\documentclass{amsart}

\usepackage{amsmath,amssymb, amsbsy}
\usepackage{color,psfrag}
\usepackage{enumerate}
\usepackage[dvips]{graphicx}
\usepackage[a4paper,width={15cm},left=3cm,bottom=2.5cm, top=2.5cm,includeheadfoot]{geometry}
\usepackage[utf8]{inputenc}

\renewcommand{\a }{\alpha }
\renewcommand{\b }{\beta }
\renewcommand{\d}{\delta }

\newcommand{\e }{\varepsilon }
\newcommand{\g }{\gamma}

\renewcommand{\l }{\lambda }
\renewcommand{\L }{\Lambda }

\newcommand{\n }{\nabla }
\newcommand{\vp }{\varphi }
\renewcommand{\phi}{\varphi}

\renewcommand{\t }{\tau }

\renewcommand{\th }{\theta }

\renewcommand{\O }{\Omega }

\newcommand{\ov}{\overline}

\newcommand{\be}{\begin{equation}}
\newcommand{\ee}{\end{equation}}
\newenvironment{pf}{\noindent{\sc Proof}.\enspace}{\rule{2mm}{2mm}\medskip}
\newenvironment{pfn}[1]{\noindent{\sc Proof of
    {#1}.\enspace}}{\hfill\qed\medskip}
\newcommand{\R}{\mathbb{R}}
\newcommand{\N}{\mathbb{N}}
\renewcommand{\S}{\mathbb{S}}

\newcommand{\de}{\partial}

\renewcommand{\k}{\kappa}

\newcommand{\calD }{\mathcal{D}}

\newcommand{\calN}{{\mathcal N}}

\newcommand{\Ds}{(-\Delta)^s}

\newcommand{\SN}{{\mathbb S}^{N-1}}
\newcommand{\weakly}{\rightharpoonup}

\newcommand{\dive }{\mathop{\rm div}}

\renewcommand{\geq }{\geqslant}

\renewcommand{\leq }{\leqslant}

\newtheorem{Theorem}{Theorem}[section]
\newtheorem{Corollary}[Theorem]{Corollary}
\newtheorem{Lemma}[Theorem]{Lemma}
\newtheorem{Proposition}[Theorem]{Proposition}
\theoremstyle{definition} 

\newtheorem{remark}[Theorem]{Remark}
\begin{document}

\title[ Fractional elliptic equations]{Unique continuation property and  local asymptotics of solutions to fractional
  elliptic equations}

 \author[Mouhamed Moustapha Fall \and Veronica Felli]{Mouhamed Moustapha Fall \and Veronica Felli }
 \address{\hbox{\parbox{5.7in}{\medskip\noindent
  M.M. Fall\\
          African Institute for Mathematical Sciences (A.I.M.S.) of Senegal,\\
         KM 2, Route de Joal, AIMS-Senegal\\
         B.P. 1418.
         Mbour, Sénégal. \\[2pt]
         {\em{E-mail address: }}{\tt mouhamed.m.fall@aims-senegal.org.}\\[5pt]
 V. Felli\\
 Universit\`a di Milano
         Bicocca,\\
         Dipartimento di Ma\-t\-ema\-ti\-ca e Applicazioni, \\
         Via Cozzi
         53, 20125 Milano, Italy. \\[2pt]
                                     \em{E-mail address: }{\tt veronica.felli@unimib.it.}}}}

\date{Revised version, June 27, 2013}

\thanks{M. M. Fall is supported by the Alexander von Humboldt foundation.
V. Felli was partially supported by the PRIN2009 grant “Critical Point Theory and Perturbative Methods for
Nonlinear Differential Equations”.
  \\
  \indent 2010 {\it Mathematics Subject Classification.} 35R11, 35B40,
    35J60, 35J75.\\
  \indent {\it Keywords.} Fractional elliptic equations,
  Caffarelli-Silvestre extension,
  Hardy inequality, Unique continuation property}

 \begin{abstract}
   \noindent Asymptotics of solutions to fractional elliptic equations
   with Hardy type potentials is studied in this paper.  By
   using an Almgren type monotonicity formula, separation of variables,
   and blow-up arguments, we describe the
   exact behavior near the singularity of solutions to linear and
   semilinear fractional  elliptic equations with a homogeneous
   singular potential related to the fractional   Hardy inequality.
   As a consequence we obtain unique continuation properties for fractional elliptic equations.
 \end{abstract}

\maketitle

\bigskip
\section{Introduction }\label{intro}

The purpose of the present paper is to describe the asymptotic behavior of
solutions to the following class of  fractional elliptic semilinear equations with
singular homogeneous potentials
\begin{equation}\label{eq:frac_eq}
(-\Delta)^su(x)-\frac{\lambda}{|x|^{2s}}u(x)=h(x)u(x)+f(x,u(x)),\quad\text{in }\Omega,
\end{equation}
where $u\in\calD^{s,2}(\R^N)$ (see definition below) and
$\Omega\subset \R^N$ is a bounded domain containing the origin,
\begin{align}
\label{eq:ipo1}&N>2s,\quad
s\in(0,1),\quad
\lambda<\Lambda_{N,s}:=2^{2s}\frac{\Gamma^2\big(\frac{N+2s}{4}\big)}{\Gamma^2\big(\frac{N-2s}{4}\big)},\\[5pt]
\label{eq:ipoh}
&h\in C^1(\Omega\setminus\{0\}),
\quad |h(x)|+|x\cdot\nabla h(x)|\leq C_h|x|^{-2s+\e} \text{ as } |x|\rightarrow 0,
\\[5pt]
\label{eq:ipof}
&\left\{\!\!
\begin{array}{l}
f\in C^1(\Omega\times \R),
\quad t\mapsto F(x,t)\in C^{1}(\Omega\times\R), \\[5pt]
|f(x,t)t|+|f'_t(x,t)t^2|+|\nabla_x F(x,t)\cdot x|\leq C_f\,|t|^{p}
\text{ for a.e. $x\in\Omega$ and all $t\in\R$},
\end{array}
\right.
\end{align}
where $2<p\leq 2^*\!(s)=\frac{2N}{N-2s}$,
$F(x,t)=\int_0^t f(x,r)\,dr$, $C_f,C_h,\e>0$ are positive
constants independent of $x\in\Omega$ and $t\in\R$, $\nabla_x F$
denotes the gradient of $F$ with respect to the $x$ variable, and
$f'_t(x,t)=\frac{\partial f}{\partial t}(x,t)$.

We recall that for any $\varphi\in C^\infty_{\rm c}(\R^N)$ and $s\in(0,1)$, the
fractional Laplacian $(-\Delta)^s\varphi$ is defined as
\begin{equation}\label{eq:2}
  (-\Delta)^s\varphi(x)=C(N,s)\mathop{\rm
    P.V.}\int_{\R^N}\frac{\varphi(x)-\varphi(y)}{|x-y|^{N+2s}}dy
=C(N,s)\lim_{\rho\to 0^+}\int_{|x-y|>\rho}\frac{\varphi(x)-\varphi(y)}{|x-y|^{N+2s}}dy
\end{equation}
where P.V. indicates  that the integral is meant in the principal
value sense
and
$$
C(N,s)=\pi^{-\frac
  N2}2^{2s}\frac{\Gamma\big(\frac{N+2s}{2}\big)}{\Gamma(2-s)}s(1-s).
$$
The Dirichlet form associated to $\Ds$ on $C^\infty_{\rm c}(\R^N)$ is given by
\begin{equation}\label{eq:3}
  (u,v)_{\mathcal
    D^{s,2}(\R^N)}=\frac{C(N,s)}2\int_{\R^{2N}}\frac{(u(x)-u(y))(v(x)-v(y))}{|x-y|^{N+2s}}\,dx\,dy
  =\int_{\R^N}
  |\xi|^{2s}\ov{\widehat v(\xi) }\widehat u(\xi)\,d\xi,
\end{equation}
where $\widehat{u}$ denotes the unitary Fourier transform of $u$. It defines a scalar product thanks to \eqref{eq:sobolev} or \eqref{eq:frac_hardy} below.
From now on we  define $\mathcal D^{s,2}(\R^N)$ as the completion
on $C^\infty_{\rm c}(\R^N)$ with respect to the norm induced by the
scalar product \eqref{eq:3}.

By a weak solution to \eqref{eq:frac_eq} we mean a function
$u\in \mathcal D^{s,2}(\R^N)$ such that
\begin{equation}\label{eq:1}
(u,\varphi)_{\mathcal D^{s,2}(\R^N)}
=\int_\Omega \bigg(\frac{\lambda}{|x|^{2s}}u(x)+h(x)u(x)+f(x,u(x))\bigg)\varphi(x)\,dx,
\text{ for all }\varphi\in C^\infty_c(\Omega).
\end{equation}

We notice that the right hand side of \eqref{eq:1} is well defined in
view of assumptions \eqref{eq:ipoh}--\eqref{eq:ipof}, the
Hardy-Littlewood-Sobolev inequality
\begin{equation}\label{eq:sobolev}
S_{N,s}  \|u\|_{L^{2^*\!(s)}(\R^N)}^2\leq\|u\|^2_{\mathcal D^{s,2}(\R^N)},
\end{equation}
and the following Hardy inequality, due to Herbst in \cite{herbst} (see also \cite{yafaev}),
\begin{equation}\label{eq:frac_hardy}
  \Lambda_{N,s}\int_{\R^N}\frac{u^2(x)}{|x|^{2s}}\,dx\leq \int_{\R^N}
  |\xi|^{2s}|\widehat u(\xi)|^2\,d\xi =\|u\|_{\mathcal D^{s,2}(\R^N)}^2
  ,\quad\text{for all }u\in \mathcal D^{s,2}(\R^N).
\end{equation}
It should also be remarked that, in \eqref{eq:1}, we allow $p=2^*(s)$  and that $u$ is not prescribed outside $\O$.

One of the aim of this paper  is to  give the precise behaviour of a solution $u$ to \eqref{eq:frac_eq}.
The rate and the shape of $u$ are given by the
the eigenvalues and the eigenfunctions of
 the following eigenvalue problem
\begin{align}\label{eq:4}
  \begin{cases}
    -\dive\nolimits_{{\mathbb S}^{N}}(\theta_1^{1-2s}\nabla_{{\mathbb
        S}^{N}}\psi)=\mu\,
    \theta_1^{1-2s}\psi, &\text{in }{\mathbb S}^{N}_+,\\[5pt]
-\lim_{\theta_1\to 0^+} \theta_1^{1-2s}\nabla_{{\mathbb
    S}^{N}}\psi\cdot {\mathbf
  e}_1=\kappa_s\lambda\psi,&\text{on }\partial {\mathbb S}^{N}_+,
  \end{cases}
\end{align}
where
\begin{equation}\label{eq:kappa_s}
\kappa_s=\frac{\Gamma(1-s)}{2^{2s-1}\Gamma(s)}
\end{equation}
and
\begin{align*}
  &{\mathbb S}^{N}_+=\{(\theta_1,\theta_2,\dots, \theta_{N+1})\in
  {\mathbb S}^{N}:\theta_1>0\}=\left\{\tfrac{z}{|z|}:z\in \R^{N+1},\ z\cdot {\mathbf e}_1>0\right\},
\end{align*}
with ${\mathbf e}_1=(1,0,\dots,0)$; we refer to section
\ref{sec:separ-vari-extens} for a variational formulation of \eqref{eq:4}. From classical
spectral theory (see section \ref{sec:separ-vari-extens} for the details)
problem \eqref{eq:4} admits a diverging sequence of real eigenvalues
with finite multiplicity
\begin{equation*}
  \mu_1(\lambda)\leq\mu_2(\lambda)\leq\cdots\leq\mu_k(\lambda)\leq\cdots
\end{equation*}
Moreover $\mu_1(\lambda)>- \big(\frac{N-2s}2\big)^2$, see Lemma
\ref{l:hardysphere}. We notice that if $\lambda=0$ then $\mu_k(0)\geq
0$ for all $k$.

Our first main result is the following theorem.
\begin{Theorem} \label{t:asym-frac}
Let $u\in \mathcal D^{s,2}(\R^N)$ be a nontrivial
solution to \eqref{eq:frac_eq} in a bounded domain $\Omega\subset
\R^N$ containing the origin as in \eqref{eq:1} with $s,\l,h$ and $f$ satisfying assumptions \eqref{eq:ipo1},
\eqref{eq:ipoh} and \eqref{eq:ipof}. Then there exists an eigenvalue $\mu_{k_0}(\l)$
of \eqref{eq:4} and an eigenfunction $\psi$ associated to $\mu_{k_0}(\l)$ such that
\be\label{eq:limint}
\t^{-\frac{2s-N}{2}-\sqrt{ \left(\frac{2s-N}{2}  \right)^2
    +\mu_{k_0}(\l)  }}u(\tau x)\to
|x|^{-\frac{N-2s}{2}+\sqrt{ \left(\frac{2s-N}{2}  \right)^2
    +\mu_{k_0}(\l)  }}
\psi\Big(0,\frac{x}{|x|}\Big)\quad\textrm{as } \t\to0^+,
\ee
in  $C^{1,\alpha}_{\rm
    loc}(B_1'\setminus\{0\})$ for some $\a\in(0,1)$,
where $B_1':=\{x\in \R^N:|x|<1\}$,
and, in particular,
\be\label{eq:limpsi}
\t^{-\frac{2s-N}{2}-\sqrt{ \left(\frac{2s-N}{2}  \right)^2 +\mu_{k_0}(\l)  }}u\left(\tau\theta'\right)\to \psi\left(0,\theta'\right)\quad\textrm{in } C^{1,\a}(\S^{N-1})\quad\textrm{as } \t\to0^+,
\ee
where $\S^{N-1}=\de \S^N_+$.
\end{Theorem}
We should mention that the above result is stated in such form for the
sake of simplicity. In fact, for the eigenfunction $\psi$ in
\eqref{eq:limpsi}, we obtain precisely its components in any basis of
the eigenspace corresponding to $\mu_{k_0}(\l) $, see Theorem
\ref{t:asym} (below).  We also remark that $\mu_{1}(\l)<0 $ for
$\l> 0$ and it is determined implicitly by the usual Gamma
function, see Proposition \ref{p:mu1la}.

Analogous results were obtained in \cite{FFT, FFT2,
  FFT3} for corresponding equations involving the Laplacian (i.e. in
the case $s=1$) and Hardy-type potentials for different kinds of
problems: in \cite{FFT, FFT3} for Schr\"odinger equations with
electromagnetic potentials and in \cite{FFT2} for Schr\"odinger
equations with inverse square many-particle potentials.

As a particular case of Theorem \ref{t:asym-frac}, if $\lambda=0$ we
obtain that the convergence stated in \ref{eq:limint} holds in  $C^{1,\alpha}(B_1')$.
\begin{Corollary} \label{t:lambda0-asym-frac}
Let $u\in \mathcal D^{s,2}(\R^N)$ be a nontrivial weak
solution to
\begin{equation}\label{eq:frac_eq_unspm}
(-\Delta)^su(x)=h(x)u(x)+f(x,u(x)),\quad\text{in }\Omega,
\end{equation}
in a bounded domain $\Omega\subset \R^N$ with $s\in(0,1)$, $h\in
C^1(\Omega)$, and $f$ satisfying \eqref{eq:ipof}. Then, for every
$x_0\in\Omega$,  there exists an eigenvalue $\mu_{k_0}=\mu_{k_0}(0)$
of problem \eqref{eq:4} with $\lambda=0$ and an eigenfunction $\psi$
associated to $\mu_{k_0}$ such that
\begin{multline}\label{eq:limint0lambda}
  \t^{-\frac{2s-N}{2}-\sqrt{ \left(\frac{2s-N}{2} \right)^2
      +\mu_{k_0}  }}u(x_0+\tau (x-x_0))\\
  \to |x-x_0|^{-\frac{N-2s}{2}+\sqrt{ \left(\frac{2s-N}{2} \right)^2
      +\mu_{k_0} }}
  \psi\Big(0,\frac{x-x_0}{|x-x_0|}\Big)\quad\textrm{as } \t\to0^+,
\end{multline}
in  $C^{1,\alpha}(\{x\in\R^N:x-x_0\in B_1'\})$.
\end{Corollary}

Our next result contains   the so called  {\em{strong unique continuation
  property}} which is a direct consequence of Theorem~\ref{t:asym-frac}.

\begin{Theorem}\label{t:sun}
  Suppose that all the assumptions of Theorem \ref{t:asym-frac} hold true.
Let $u$ be a solution to \eqref{eq:frac_eq} in a bounded domain
$\Omega\subset \R^N$ containing the origin. If
$u(x)=o(|x|^n)=o(1)|x|^n$ as $|x|\to 0$ for all $n\in
  \N$, then $u\equiv 0$ in $\Omega$.
\end{Theorem}

From Theorem \ref{t:lambda0-asym-frac} a  unique continuation
principle related to sets of positive measures follows.
\begin{Theorem}\label{t:unspm}
Let $u\in \mathcal D^{s,2}(\R^N)$ be a weak
solution to \eqref{eq:frac_eq_unspm}
in a bounded domain $\Omega\subset
\R^N$ with $s\in(0,1)$, $h\in C^1(\Omega)$, and $f$ satisfying
\eqref{eq:ipof}.
If $u\equiv 0$ on a set $E\subset \Omega$ of positive measure, then $u\equiv 0$ in $\Omega$.
\end{Theorem}

An interesting application of Theorem \ref{t:unspm} is that the
nodal sets of eigenfunctions for the fractional laplacian operator
have zero Lebesgue measure. We should point out that this was a key
assumption in \cite{FSV}, where the authors studied  the existence
of weak solutions to some non-local equations.

Recent research in the field of second order elliptic equations has
devoted a great attention to the problem of unique continuation
property in the presence of singular lower order terms, see
e.g. \cite{FGL,JK,kurata}.  Two different kinds of approach have been
developed to treat unique continuation: a first one is due to Carleman
\cite{carleman} and is based on weighted priori inequalities, whereas
a second one is due to Garofalo and Lin \cite{GL} and is based on
local doubling properties proved by Almgren monotonicity formula. In
the present paper we will follow the latter approach. Furthermore,
in the spirit  of \cite{FFT, FFT2, FFT3}, the combination of
monotonicity methods with blow-up analysis will enable us to prove not
only unique continuation but also the
precise asymptotics of solutions stated in Theorem \ref{t:asym-frac}.

As far as  unique continuation from sets of positive measures is
concerned, we mention \cite{deFG} where it was proved for second
order elliptic operators by
 combining strong unique continuation
property with the De Giorgi inequality. Since the validity of a the
De Giorgi type inequality for the fractional problem (or even its
extension, see \eqref{eq:Harmext}) seems to be hard to prove, we
will base the proof of Theorem \ref{t:unspm} directly on the
asymptotic of solutions proved in Theorem \ref{t:asym-frac}.

To explain our argument of proving the asymptotic behavior,  let us
we write \eqref{eq:1} as
\begin{equation}\label{eq:50}
\Ds u=G(x,u)\quad \textrm{in } \O.
\end{equation}
The proof of our results is based on the study of the \textit{Almgren
  frequency function} at the origin 0: ``ratio of the local energy
over mass near the origin''.  Due to the non-locality of the Dirichlet
form associated to $\Ds$, it is not clear how to set up an Almgren's
type frequency function using this energy as in the local case
$s=1$. A way out for this difficulty is to use the
Caffarelli-Silvestre extension \cite{CS} which can be seen as a local
version of \eqref{eq:50}.  The Caffarelli-Silvestre extension of a
solution $u$ to \eqref{eq:50} is a function $w$ defined on
\[
\R^{N+1}_+=\{z=(t,x):t\in(0,+\infty),\ x\in\R^N\}
\]
satisfying $w =u$ on $\O$ and solving in some weak sense (see
Section~\ref{sec:some-preliminaries} for more details) the boundary
value problem
\begin{equation}\label{eq:Harmext}
\begin{cases}
\dive(t^{1-2s}\n w)=0,&\textrm{ in }\R^{N+1}_+,\\
-  \lim \limits_{t \to 0^+}t^{1-2s}\,\frac{\de w}{\de t}= \k_sG(x,w), &\textrm{ on } \O.
\end{cases}
\end{equation}
Here and in the following, we write $z=(t,x) \in \R^{N+1}_+$ with
$x \in \R^N$ and $t>0$, and we identify $\R^N$ with $\partial \R^{N+1}_+$, so that
$\Omega$ is contained in $\partial \R^{N+1}_+$.

We then consider the {Almgren's frequency} function
\begin{equation*}
\calN(r)=\frac{D(r)}{H(r)}
\end{equation*}
where
\begin{equation*}
D(r)=\frac{1}{r^{N-2s}}\bigg[\int_{B_r^+}t^{1-2s}|\n
w|^2\,dt\,dx-\kappa_s\int_{B_r'}G(x,w)w\, dx\bigg],\quad
H(r)=\frac{1}{r^{N+1-2s}}\int_{S_r^+}t^{1-2s}w^2 \, dS,
\end{equation*}
being
\begin{align*}
&B_r^+=\{z=(t,x)\in\R^{N+1}_+\,:\,  |z|<r\} ,\quad
B_r':=\{x\in \R^N:|x|<r\},\\
&S_r^+=\{z=(t,x)\in\R^{N+1}_+ \,:\, |z|=r\},
\end{align*}
and $dS$ denoting the volume element on $N$-dimensional spheres. We
mention that an Almgren's frequency function for degenerate elliptic
equations of the form \eqref{eq:Harmext} was first formulated in
\cite[Section 6]{CS}.

As a first but nontrivial step, we prove that $\lim _{r\to
  0}\calN(r):=\g$ exists and it is finite, see Lemma~\ref{gamma}.
Next,  we make a blow-up analysis by zooming around the origin the
solution $w$ normalized also by $\sqrt{H}$.  More precisely, setting
$w_\t(z)=\frac{w(\t z)}{\sqrt{H(\t)}}$, we have that $w_\t$ converges, (in some
H\"{o}lder and Sobolev spaces) to $\widetilde{w}$ solving the limiting
equation
\begin{equation*}
\begin{cases}
\dive(t^{1-2s}\n \widetilde{w})=0, &\textrm{in }B_1^+,\\
- \lim \limits_{t \to 0^+}t^{1-2s}\,\frac{\de \widetilde{w}}{\de t}=
\frac{\k_s \l}{|z|^{2s}}\widetilde{w},& \textrm{on } B_1'.
\end{cases}
\end{equation*}
To obtain this, the fact that $h$ is negligible with respect  to
the Hardy potential and $f$ is at most critical with respect to
the Sobolev exponent (see assumptions \eqref{eq:ipoh} and
\eqref{eq:ipof}) plays a crucial role;
we refer to Lemma \ref{l:blowup} for more details.

The main point is that the Almgren's frequency for $\widetilde{w}$ is
$\widetilde{\calN}(r)=\lim_{\t\to 0}\calN(\tau r)=\g$,
i.e. $\widetilde{\calN}$ is constant; hence
$\widetilde{w}$ and $\n\widetilde{w}\cdot\frac{z}{|z|}$ are proportional on $L^2(S^+_r;t^{1-2s})$.
As a consequence we obtain $\widetilde{w}(z)=\vp(|z|)\psi\big( \frac{z}{|z|}\big)$.
By separating variables in polar coordinates, we obtain that $\psi$ is an eigenfunction
of \eqref{eq:4} for
some eigenvalue $\mu_{k_0}(\l)$.
By the method of  variation of constants and the fact that $\widetilde{w}$
has finite energy near the origin, we then prove that $\varphi(r)$ is
proportional to $r^{\g}$ and that $\gamma=\frac{2s-N}{2}+\sqrt{ \left(\frac{2s-N}{2}  \right)^{\!2} +\mu_{k_0}(\l)  }$.
We finally complete the proof by showing that $\lim_{r\to0} r^{-2\g}H(r)>0$, see Lemma \ref{l:limitepositivo}.\\

We should mention that in the recent literature, a great attention has
been addressed to nonlocal fractional diffusion and many papers have
been devoted to the study of existence, non-existence, regularity and
qualitative properties of solutions to elliptic equations associated
to fractional Laplace type operators, see
e.g. \cite{BCdP,CaSi,CS,DNPV,fall-frac,FW,silvestre,SV} and references
therein. In particular, semilinear fractional elliptic equations
involving the Hardy potential were treated in \cite{fall-frac}, where
some existence and nonexistence results were obtained from lower
bounds of positive solutions.  Adapting the ideas in \cite{FFT} in the
nonlocal case of the present paper, we provide the behavior (therefore
regularity) of solutions at the singularity.  In particular, the
asymptotics we have here show that the lower bound of positive
solutions proved in \cite{fall-frac} is sharp.

\section{Preliminaries and notations}\label{sec:some-preliminaries}

\medskip
\noindent
{\bf Notation. } We list below some notation used throughout the
paper.\par
\begin{itemize}
\item[-] ${\mathbb S}^{N}=\{z\in \R^{N+1}:|z|=1\}$ is the unit
  $N$-dimensional sphere.
\item[-] ${\mathbb S}^{N}_+=\{(\theta_1,\theta_2,\dots, \theta_{N+1})\in
  {\mathbb S}^{N}:\theta_1>0\}:={\mathbb S}^{N}\cap\R^{N+1}_+$.
\item[-] $dS$ denotes the volume
element on $N$-dimensional spheres.
\item[-] $dS'$ denotes the volume
element on $(N-1)$-dimensional spheres.
\end{itemize}
Let $\mathcal D^{1,2}(\R^{N+1}_+;t^{1-2s})$ be the completion of
$C^\infty_{\rm c}(\overline{\R^{N+1}_+})$ with respect to the norm
$$
\|w\|_{\mathcal D^{1,2}(\R^{N+1}_+;t^{1-2s})}=\bigg(
\int_{\R^{N+1}_+}t^{1-2s}|\nabla w(t,x)|^2dt\,dx
\bigg)^{\!\!1/2}.
$$
We recall that there exists a well defined continuous trace map
$\mathop{\rm Tr}:\mathcal D^{1,2}(\R^{N+1}_+;t^{1-2s})\to \mathcal
\calD^{s,2}(\R^N)$, see e.g. \cite{BCdP}.

For every $u\in \mathcal D^{s,2}(\R^N)$, let  $\mathcal H(u) \in
\mathcal D^{1,2}(\R^{N+1}_+;t^{1-2s})$ be the
unique solution to the minimization problem
\begin{multline}\label{eq:5}
  \int_{\R^{N+1}_+}t^{1-2s}|\nabla \mathcal
  H(u)|^2\,dt\,dx\\
  =\min\bigg\{ \int_{\R^{N+1}_+}t^{1-2s}|\nabla v|^2\,dt\,dx:v\in
  \mathcal D^{1,2}(\R^{N+1}_+;t^{1-2s}), \ \mathop{\rm
    Tr}(v)=u\bigg\}.
\end{multline}
By Caffarelli and Silvestre \cite{CS} we have that
\begin{equation}\label{eq:6}
\int_{\R^{N+1}_+}t^{1-2s}\nabla \mathcal
H(u)\cdot\nabla\widetilde\varphi\,dt\,dx=
\kappa_s(u,\mathop{\rm Tr}\widetilde\varphi)_{\mathcal
  D^{s,2}(\R^N)}\quad\text{for all }\widetilde\varphi\in  \mathcal D^{1,2}(\R^{N+1}_+;t^{1-2s}),
\end{equation}
where $\kappa_s$ is defined in \eqref{eq:kappa_s}.

We notice that combining \eqref{eq:frac_hardy}, \eqref{eq:5},
\eqref{eq:6} and \eqref{eq:sobolev}, we obtain the following Hardy-trace inequality
 \begin{align}\label{eq:half_space_hardy}
   \kappa_s\Lambda_{N,s}\int_{\R^N}\frac{(\mathop{\rm
       Tr}v)^2(x)}{|x|^{2s}}\,dx&\leq \int_{\R^{N+1}_+}t^{1-2s}|\nabla
   v|^2\,dt\,dx,\quad\text{for all }v\in{\mathcal
     D}^{1,2}({\R^{N+1}_+};t^{1-2s})
 \end{align}
 and also the Sobolev-trace inequality
  \begin{align}\label{eq:half_space_Sobolev}
   \kappa_sS_{N,s}  \|\mathop{\rm
       Tr}v\|_{L^{2^*\!(s)}(\R^N)}^2
        &\leq \int_{\R^{N+1}_+}t^{1-2s}|\nabla
   v|^2\,dt\,dx,\quad\text{for all }v\in{\mathcal
     D}^{1,2}({\R^{N+1}_+};t^{1-2s}).
 \end{align}

\subsection{Separation of variables in the extension operator}\label{sec:separ-vari-extens}
For every $R>0$, we  define the space
$H^1(B_R^+;t^{1-2s})$ as the  completion of
$C^\infty(\overline{B_R^+})$ with respect to the norm
$$
\|w\|_{H^1(B_R^+;t^{1-2s})}=\bigg(
\int_{B_R^+}t^{1-2s}\Big(|\nabla w(t,x)|^2+w^2(t,x)\Big)
dt\,dx
\bigg)^{\!\!1/2}.
$$

By direct calculations, we can obtain the following lemma concerning
separation of variables in the extension operator.
\begin{Lemma}\label{l:separation}
  If  $v\in H^1(B_R^+;t^{1-2s})$ is such that
$v(z)=f(r)\psi(\theta)$ for a.e. $z=(t,x)\in \R_+^{N+1}$, with
$r=|z|<R$ and $\theta=\frac{z}{|z|}\in {\mathbb S}^{N}$, then
$$
 \dive(t^{1-2s}\nabla
 v(z))=\frac1{r^N}\big(r^{N+1-2s}f'\big)'\theta_1^{1-2s}\psi(\theta)+
r^{-1-2s}f(r)\dive\nolimits_{{\mathbb S}^{N}}(\theta_1^{1-2s}\nabla_{{\mathbb S}^{N}}\psi(\theta))
$$
in the distributional sense, where $\theta_1=\frac{t}{r}=\theta\cdot
{\mathbf e}_1$ with ${\mathbf e}_1=(1,0,\dots,0)$,
$\dive\nolimits_{{\mathbb S}^{N}}$ (respectively $\nabla_{{\mathbb S}^{N}}$) denotes the Riemannian
 divergence (respectively gradient) on the unit sphere ${\mathbb S}^{N}$ endowed
 with the standard metric.
\end{Lemma}

\noindent
Let us define $H^{1}({\mathbb S}^{N}_+;\theta_1^{1-2s})$ as the completion of
$C^\infty(\overline{{\mathbb S}^{N}_+})$ with respect to the norm
$$
\|\psi\|_{H^1({\mathbb S}^{N}_+;\theta_1^{1-2s})}=\bigg(
\int_{{\mathbb S}^{N}_+}\theta_1^{1-2s}\big(|\nabla_{{\mathbb
        S}^{N}}\psi(\theta)|^2+\psi^2(\theta)\big)dS
\bigg)^{\!\!1/2}.
$$
We also denote
$$
L^2({\mathbb
  S}^{N}_+;\theta_1^{1-2s}):=\Big\{\psi:{\mathbb S}_+^{N}\to\R\text{
  measurable such that }{\textstyle{\int}}_{{\mathbb
    S}^{N}_+}\theta_1^{1-2s}\psi^2(\theta)\,dS<+\infty\Big\}.
$$
The following trace inequality on the unit half-sphere ${\mathbb S}_+^{N}$ holds.
\begin{Lemma}\label{l:hardysphere}
  There exists a well defined continuous trace operator
$$
 H^1({\mathbb S}^{N}_+;\theta_1^{1-2s})\to L^2(\partial{\mathbb
    S}^{N}_+)=L^2({\mathbb S}^{N-1}).
$$
 Moreover for every $\psi\in
  H^1({\mathbb S}^{N}_+;\theta_1^{1-2s})$
\begin{align*}
 \kappa_s\Lambda_{N,s}\bigg( \int_{{\mathbb
      S}^{N-1}}| \psi(\theta')|^2\,dS'\bigg)\leq \Big(\frac{N-2s}2\Big)^2  \int_{{\mathbb
      S}^{N}_+}\theta_1^{1-2s}|\psi(\theta)|^2\,dS+
\int_{{\mathbb
      S}^{N}_+}\theta_1^{1-2s}|\nabla_{{\mathbb
      S}^{N}} \psi(\theta)|^2\,dS.
\end{align*}
where $dS'$ denotes the volume
element on the sphere ${\mathbb S}^{N-1}=\partial {\mathbb
 S}^{N}_+=\{(\theta_1,\theta')\in {\mathbb S}^{N}_+:\theta_1=0\}$.
\end{Lemma}
\begin{pf}
Let $\psi\in C^\infty(\overline{{\mathbb S}^{N}_+})$ and $f\in C^\infty_{\rm
c}(0,+\infty)$ with $f\neq0$. Rewriting \eqref{eq:half_space_hardy} for
$v(z)=f(r)\psi(\theta)$, $r=|z|$, $\theta=\frac z{|z|}$, we
obtain that
\begin{align*}
  \kappa_s\Lambda_{N,s}&\bigg(\int_0^{+\infty
  }r^{N-1-2s}\,f^2(r)\,dr\bigg)\bigg( \int_{{\mathbb
      S}^{N-1}}|\psi(0,\theta')|^2\,dS'\bigg)\\
&\leq \bigg(\int_0^{+\infty
  }r^{N+1-2s}|f'(r)|^2\,dr\bigg)\bigg( \int_{{\mathbb
      S}^{N}_+}\theta_1^{1-2s}|\psi(\theta)|^2\,dS\bigg)\\
&\quad+\bigg(\int_0^{+\infty
  }r^{N-1-2s}f^2(r)\,dr\bigg)\bigg( \int_{{\mathbb
      S}^{N}_+}\theta_1^{1-2s}|\nabla_{{\mathbb
      S}^{N}_+} \psi(\theta)|^2\,dS\bigg),
\end{align*}
and hence, by optimality of the classical Hardy constant, see
\cite[Theorem 330]{HLP},
\begin{align*}
 \kappa_s\Lambda_{N,s}&\bigg( \int_{{\mathbb
      S}^{N-1}}|\psi(0,\theta')|^2\,dS'\bigg)
\\
&\leq \bigg( \int_{{\mathbb
      S}^{N}_+}\theta_1^{1-2s}|\psi(\theta)|^2\,dS\bigg)
\inf_{f\in C^\infty_{\rm c}(0,+\infty)}
\frac{\int_0^{+\infty
  }r^{N+1-2s}|f'(r)|^2\,dr}{\int_0^{+\infty
  }r^{N-1-2s}\,f^2(r)\,dr}
+ \int_{{\mathbb
      S}^{N}_+}\theta_1^{1-2s}|\nabla_{{\mathbb
      S}^{N}_+} \psi(\theta)|^2\,dS\\
&=\Big(\frac{N-2s}2\Big)^2  \int_{{\mathbb
      S}^{N}_+}\theta_1^{1-2s}|\psi(\theta)|^2\,dS+
\int_{{\mathbb
      S}^{N}_+}\theta_1^{1-2s}|\nabla_{{\mathbb
      S}^{N}_+} \psi(\theta)|^2\,dS.
\end{align*}
By density of $C^\infty(\overline{{\mathbb S}^{N}_+})$ in
$H^1({\mathbb S}^{N}_+;\theta_1^{1-2s})$, we obtain the conclusion.
\end{pf}

In view of Lemma \ref{l:separation}, in order to construct an
orthonormal basis of $L^2({\mathbb S}^{N}_+;\theta_1^{1-2s})$ for
expanding solutions to \eqref{eq:Harmext} in Fourier series, we are
naturally lead to consider the eigenvalue problem \eqref{eq:4}, which
admits the following variational formulation: we say that $\mu\in\R$
is an eigenvalue of problem \eqref{eq:4} if there exists $\psi\in H^1({\mathbb
  S}^{N}_+;\theta_1^{1-2s})\setminus\{0\}$ (called eigenfunction) such
that
$$
Q(\psi,\upsilon)=\mu\int_{{\mathbb
    S}^{N}_+}\theta_1^{1-2s}\psi(\theta)\upsilon(\theta)\,dS,\quad\text{for
all }\upsilon \in H^1({\mathbb S}^{N}_+;\theta_1^{1-2s}),
$$
where
\begin{align*}
&Q:H^1({\mathbb S}^{N}_+;\theta_1^{1-2s})\times H^1({\mathbb
  S}^{N}_+;\theta_1^{1-2s})\to\R,\\
&Q(\psi,\upsilon)=\int_{{\mathbb S}^{N}_+}\theta_1^{1-2s}\nabla_{{\mathbb
        S}^{N}}\psi(\theta)\cdot\nabla\upsilon(\theta)\,dS-\lambda
    \kappa_s \int_{{\mathbb
      S}^{N-1}}\mathcal T\psi(\theta') \mathcal T\upsilon(\theta')\,dS' .
\end{align*}
By Lemma \ref{l:hardysphere} the bilinear form $Q$ is continuous and
weakly coercive on $H^1({\mathbb S}^{N}_+;\theta_1^{1-2s})$. Moreover the
belonging of the weight $t^{1-2s}$ to the second Muckenhoupt
class ensures  that the embedding $H^1({\mathbb
  S}^{N}_+;\theta_1^{1-2s})\hookrightarrow \hookrightarrow
L^2({\mathbb S}^{N}_+;\theta_1^{1-2s})$ is compact (see \cite{FKS} for
weighted embeddings with  Muckenhoupt $A_2$ weights). Then, from classical
spectral theory (see e.g. \cite[Theorem 6.16]{salsa}),
problem \eqref{eq:4} admits a diverging sequence of real eigenvalues
with finite multiplicity
  $\mu_1(\lambda)\leq\mu_2(\lambda)\leq\cdots\leq\mu_k(\lambda)\leq\cdots$ the
  first of which admits the variational  characterization
\begin{equation}\label{firsteig}
  \mu_1(\lambda)=\min_{\psi\in H^1({\mathbb
  S}^{N}_+;\theta_1^{1-2s}) \setminus\{0\}}\frac{Q(\psi,\psi)}{\int_{{\mathbb
    S}^{N}_+}\theta_1^{1-2s}\psi^2(\theta)\,dS}.
\end{equation}
Furthermore, in view of Lemma \ref{l:hardysphere}, we have that
\begin{equation}\label{firsteig_strict_in}
\mu_1(\lambda)>- \Big(\frac{N-2s}2\Big)^2.
\end{equation}
To each $k\geq 1$, we
associate an $L^2({\mathbb
  S}^{N}_+;\theta_1^{1-2s})$-normalized
eigenfunction $\psi_k\in H^1({\mathbb
  S}^{N}_+;\theta_1^{1-2s})\setminus\{0\}$ corresponding to
 the $k$-th eigenvalue $\mu_{k}(\lambda)$,
i.e. satisfying
\begin{equation}\label{angular}
Q(\psi_k,\upsilon)=\mu_k(\lambda)\int_{{\mathbb
    S}^{N}_+}\theta_1^{1-2s}\psi_k(\theta)\upsilon(\theta)\,dS,\quad\text{for
all }\upsilon \in H^1({\mathbb S}^{N}_+;\theta_1^{1-2s}).
\end{equation}
In the enumeration
$\mu_1(\lambda)\leq\mu_2(\lambda)\leq\cdots\leq\mu_k(\lambda)\leq
\cdots$ we repeat each eigenvalue as many times as its multiplicity;
thus exactly one eigenfunction $\psi_k$ corresponds to each index
$k\in{\mathbb{N}}$, $k\geq1$. We can choose the functions $\psi_k$ in
such a way that they form an orthonormal basis of $L^2({\mathbb
  S}^{N}_+;\theta_1^{1-2s}) $.

We can also determine $\mu_1(\lambda)$ for
$\lambda\in(0,\L_{N,s})$, where
$\Lambda_{N,s}$ is the fractional Hardy constant  defined in \eqref{eq:ipo1}.
\begin{Proposition}\label{p:mu1la}
 For every $\a\in\left(0, \frac{N-2s}{2}\right)$, we define
 $$
\l(\a)=2^{2s}  \frac{\Gamma\left(\frac{N+2s+2\a}{4}\right)}{\Gamma\left(\frac{N-2s-2\a}{4}\right)}
\frac{\Gamma\left(\frac{N+2s-2\a}{4}\right)}{\Gamma\left(\frac{N-2s+2\a}{4}\right)}.
$$
Then the mapping $\a\mapsto\l(\a)$ is continuous and decreasing. In addition we have that
$$
\mu_1(\l(\a))=\a^ 2-\left( \frac{N-2s}{2} \right)^ 2\quad \text{for
  all }\a\in   \left(0, \frac{N-2s}{2}\right).
$$
\end{Proposition}
\begin{pf}
It was proved in \cite[Lemma 3.1]{fall-frac} that, for every $\a\in\left(0, \frac{N-2s}{2}\right)$, there exists
a positive continuous function
$\Phi_{\a}:\R^{N+1}_+\to \R$  such that
\be\label{eq:v-al_def}
\begin{cases}
\dive(t^{1-2s}\n \Phi_\a)=0\quad \textrm{ in } \R^{N+1}_+\\
\Phi_\a=|x|^{ \frac{2s-N}{2} +\a}\quad  \textrm{ on } \R^{N}\setminus\{0\}\\
-t^{1-2s}\frac{\de \Phi_\a }{\de t}
=\k_s\l(\a) |{x}|^{-2s}\,\Phi_\a\quad\textrm{ on } \R^{N}\setminus\{0\},
\end{cases}
\ee
Moreover   $\Phi_\a\in H^1(B_R^+;t^{1-2s})$ for every $R>0$ and
$\Phi_\alpha$ is scale invariant, i.e.
\begin{equation*}
\Phi_\a(\tau z)= \tau^{
  \frac{2s-N}{2}+\a}\Phi_\a(z),\quad\text{for all $\tau>0$},
\end{equation*}
thus  implying that,
for all $z\in \R^{N+1}_+$,
\be\label{eq:estPsa}
c_1|z|^{ \frac{2s-N}{2}+\a} \leq \Phi_\a(z) \leq c_2 |z|^{ \frac{2s-N}{2}+\a},
\ee
for some positive constants $c_1,c_2$.
It is also known (see for instance \cite{FLS}) that the map $\a\mapsto \l(\a)$
is continuous and monotone decreasing.

We write
$$
\Phi_\a(z)=\sum_{k=1}^\infty\Phi_\a^k(|z|)\psi_k(z/|z|), \quad \Phi_\a^k(r)=\int_{\S^ {N}_+}\psi_k(\th)\Phi_\a(r\th) dS.
$$
In particular,  since  $\psi_1>0$, by \eqref{eq:estPsa} we have, for every $r\in(0,R)$,
\be\label{eq:estPsa_1}
c_1'|r|^{ \frac{2s-N}{2}+\a} \leq \Phi_\a^1(r) \leq c_2' |r|^{ \frac{2s-N}{2}+\a},
\ee
for some positive constants $c_1',c_2'$.
Using \eqref{eq:v-al_def}, we have, weakly, for every $k\geq1$ and $r\in(0,R)$,
$$
\begin{cases}
\frac1{r^N}\big(r^{N+1-2s}( \Phi_\a^k)'\big)'\theta_1^{1-2s}\psi_k(\theta)+
r^{-1-2s} \Phi_\a^k\dive\nolimits_{{\mathbb
    S}^{N}}(\theta_1^{1-2s}\nabla_{{\mathbb S}^{N}}\psi_k(\theta))
=0,\\
- \Phi_\a^k\lim_{\theta_1\to 0^+} \theta_1^{1-2s}\nabla_{{\mathbb
    S}^{N}}\psi_k(\theta)\cdot {\mathbf
  e}_1=\kappa_s\,\lambda(\a)\,\psi_k(0,\theta') \Phi_\a^k.
\end{cases}
$$
Testing the above equation with $\psi_1>0$ and using also the fact that $\Phi_\a^1>0$, we
obtain
$$
( \Phi_\a^1)''+\frac{N+1-2s}{r}( \Phi_\a^1)'-\frac{\mu_1(\l(\a))}{r^2}\, \Phi_\a^1=0
$$
and hence $\Phi_\a^1(r)$ is of the form
$$
\Phi_\a^1(r)=c_3 r^{\sigma^+}+c_4 r^{\sigma^-}
$$
for some $c_3,c_4\in\R$, where
\begin{equation*}
  \sigma^+=-\frac{N-2s}{2}+\sqrt{\bigg(\frac{N-2s}
    {2}\bigg)^{\!\!2}+\mu_{1}(\lambda(\alpha))}\quad\text{and}\quad
  \sigma^-=-\frac{N-2s}{2}-\sqrt{\bigg(\frac{N-2s}{2}
    \bigg)^{\!\!2}+\mu_1(\lambda(\a)
)}.
\end{equation*}
Since the function
 $|z|^{\sigma_{k_0}^-}\psi_1(\frac{z}{|z|})
\notin  H^1(B_1^+;t^{1-2s})$, we deduce that $c_4=0$ and thus for every $r\in(0,R)$
$$
\Phi_ \a^1(r)=c_3
r^{\sigma^+}.
$$
This together with \eqref{eq:estPsa_1} implies that
$$
\mu_1(\l(\a))=\a^ 2-\left( \frac{N-2s}{2} \right)^ 2\quad \text{for all }\a\in   \left(0, \frac{N-2s}{2}\right),
$$
as claimed.
\end{pf}

\subsection{Hardy type inequalities}\label{sec:hardy-type-ineq}
From well-known weighted embedding inequalities and the fact that the
weight $t^{1-2s}$ belongs to the second Muckenhoupt class (see
e.g. \cite{FKS}),
the embedding $H^1(B_r^+;t^{1-2s})\hookrightarrow L^2(B_r^+;t^{1-2s})$
is compact. It can be also  proved that both  the trace operators
\begin{equation}\label{eq:tracecompact}
H^1(B_r^+;t^{1-2s})\hookrightarrow \hookrightarrow L^2(S_r^+;\theta_1^{1-2s}),
\end{equation}
\begin{equation}\label{eq:tracecompactB}
H^1(B_r^+;t^{1-2s})\hookrightarrow \hookrightarrow L^2(B_r')
\end{equation}
are  well defined and compact.

\smallskip
For sake of simplicity, in the following of this paper, we will often
denote the trace of a function with the same letter as the function
itself.

\smallskip
\noindent
The following Hardy type inequality with boundary terms holds.

\begin{Lemma}\label{l:hardy_boundary}
  For all $r>0$ and $w\in H^1(B_r^+;t^{1-2s})$, the following
  inequality holds
  \begin{align*}
    \bigg(\frac{N-2s}{2}\bigg)^{\!\!2}\int_{B_r^+}t^{1-2s}\frac{w^2(z)}{|z|^2}\,dz
\leq \int_{B_r^+}t^{1-2s}\bigg(\nabla w(z)\cdot\frac{z}{|z|}\bigg)^{\!\!2}dz
+    \bigg(\frac{N-2s}{2r}\bigg)\int_{S_r^+}t^{1-2s}w^2dS.
  \end{align*}
\end{Lemma}
\begin{pf}
By scaling, it is enough to prove the stated inequality for $r=1$. Let
$$
V(z)=|z|^{\frac{2s-N}{2}},\quad z\in \R^{N+1}_+\setminus\{0\}.
$$
We notice that $V$ satisfies
\begin{equation}\label{eq:7}
-\dive(t^{1-2s}\nabla
V)=\bigg(\frac{N-2s}{2}\bigg)^{\!\!2}t^{1-2s}|z|^{-2}V(z)\quad\text{in
}\R^{N+1}_+\setminus\{0\}.
\end{equation}
Hence, letting $w\in C^\infty(\overline{B_1^+})$,  multiplying
\eqref{eq:7} with $\frac{w^2}{V}$, and integrating over
$B_1^+\setminus B_\delta^+$ with $\delta\in(0,1)$, we obtain
\begin{multline*}
 \bigg(\frac{N-2s}{2}\bigg)^{\!\!2}\int_{B_1^+\setminus B_\delta^+}t^{1-2s}\frac{w^2(z)}{|z|^2}\,dz\\
= \int_{B_1^+\setminus B_\delta^+}t^{1-2s}\nabla
V(z)\cdot\nabla\Big(\frac{w^2}{V}\Big)(z)\,dz-\int_{S_1^+}t^{1-2s}(\nabla
V\cdot\nu) \frac{w^2}{V}dS+\int_{S_\delta^+}t^{1-2s}(\nabla
V\cdot\nu)
\frac{w^2}{V}dS\\
=
-(N-2s)\int_{B_1^+\setminus B_\delta^+}t^{1-2s}\frac{w}{|z|}\Big(\nabla w\cdot \frac{z}{|z|}\Big)\,dz
-\bigg(\frac{N-2s}{2}\bigg)^{\!\!2}\int_{B_1^+\setminus B_\delta^+}t^{1-2s}\frac{w^2}{|z|^2}\,dz\\
+\frac{N-2s}{2}\int_{S_1^+}t^{1-2s}w^2dS-\frac{N-2s}{2}\frac1\delta
\int_{S_\delta^+}t^{1-2s}w^2dS
\end{multline*}
where $\nu(z)=\frac{z}{|z|}$. Since, by  Schwarz's inequality
\begin{multline}
  -(N-2s)\int_{B_1^+\setminus
    B_\delta^+}t^{1-2s}\frac{w}{|z|}\Big(\nabla w\cdot
  \frac{z}{|z|}\Big)\,dz\\
  \leq \bigg(\frac{N-2s}{2}\bigg)^{\!\!2}\int_{B_1^+\setminus
    B_\delta^+}t^{1-2s}\frac{w^2}{|z|^2}\,dz+ \int_{B_1^+\setminus
    B_\delta^+}t^{1-2s}\Big(\nabla w\cdot \frac{z}{|z|}\Big)^2\,dz
\end{multline}
and
$$
\frac1\delta
\int_{S_\delta^+}t^{1-2s}w^2dS=O(\delta^{N-2s})=o(1)\quad\text{for }\delta\to0 ^+,
$$
letting $\delta\to 0$ we obtain the stated inequality for $r=1$ and
$w\in C^\infty(\overline{B_1^+})$. The conclusion follows by density
of $C^\infty(\overline{B_1^+})$ in $H^1(B_1^+;t^{1-2s})$.
\end{pf}

\begin{Lemma}\label{l:hardy_boundary2}
For every $r>0$ and $w\in H^1(B_r^+;t^{1-2s})$, the following
  inequality holds
  \begin{align*}
\kappa_s    \Lambda_{N,s}\int_{B_r'}\frac{
  w^2}{|x|^{2s}}\,dx&\leq
\bigg(\frac{N-2s}{2r}\bigg)\int_{S_r^+}t^{1-2s}w^2 dS+
\int_{B_r^+}t^{1-2s}|\nabla w|^2dz.
  \end{align*}
\end{Lemma}
\begin{pf}
Let $w\in C^\infty(\overline{B_r^+})$. Then, passing to polar
coordinates and using Lemmas \ref{l:hardysphere} and \ref{l:hardy_boundary}, we obtain
\begin{align*}
  \kappa_s \Lambda_{N,s}&\int_{B_r'}\frac{w^2(0,x)}{|x|^{2s}}\,dx=
  \kappa_s \Lambda_{N,s}\int_0^r \rho^{N-1-2s}\bigg( \int_{{\mathbb
      S}^{N-1}}w^2(0,\rho\theta') \,dS' \bigg)d\rho\\
  &\leq \int_0^r \rho^{N-1-2s}\bigg( \Big(\frac{N-2s}2\Big)^2
  \int_{{\mathbb S}^{N}_+}\theta_1^{1-2s}|w(\rho\theta)|^2\,dS+
  \int_{{\mathbb S}^{N}_+}\theta_1^{1-2s}|\nabla_{{\mathbb
      S}^{N}} w(\rho\theta)|^2\,dS\bigg)d\rho\\
  &=\Big(\frac{N-2s}2\Big)^2\int_{B_r^+}t^{1-2s}\frac{w^2}{|z|^2}\,dz+
  \int_0^r \rho^{N-1-2s}\bigg(\int_{{\mathbb
      S}^{N}_+}\theta_1^{1-2s}|\nabla_{{\mathbb
      S}^{N}} w(\rho\theta)|^2\,dS\bigg)d\rho\\
  &\leq
  \bigg(\frac{N-2s}{2r}\bigg)\int_{S_r^+}t^{1-2s}w^2dS\\
&\hskip3cm +\int_0^r
  \rho^{N+1-2s}\bigg(\int_{{\mathbb
      S}^{N}_+}\theta_1^{1-2s}\bigg(
\frac1{\rho^2}|\nabla_{{\mathbb S}^{N}}
  w(\rho\theta)|^2+\Big|\frac{\partial w}{\partial \rho}(\rho\theta)\Big|^2\bigg)
\,dS\bigg)d\rho\\
&=\bigg(\frac{N-2s}{2r}\bigg)\int_{S_r^+}t^{1-2s}w^2dS+
\int_{B_r^+}t^{1-2s}|\nabla w|^2dz.
\end{align*}
 The conclusion follows by density
of $C^\infty(\overline{B_r^+})$ in $H^1(B_r^+;t^{1-2s})$.
\end{pf}

The following Sobolev type inequality with boundary terms holds.

\begin{Lemma}\label{l:sobolev_boundary}
  There exists $\widetilde S_{N,s}>0$ such that, for all $r>0$ and $w\in
  H^1(B_r^+;t^{1-2s})$,
  \begin{align*}
    \bigg(\int_{B_r'}| w |^{2^*(s)}\,dx\bigg)^{\!\!\frac2{2^*(s)}}
\leq
\widetilde S_{N,S}\bigg[
\frac{N-2s}{2r}\int_{S_r^+}t^{1-2s}w^2dS+
\int_{B_r^+}t^{1-2s}|\nabla w|^2dz\bigg].
  \end{align*}
\end{Lemma}
\begin{pf}
By scaling, it is enough to prove the statement  for $r=1$.
Let $w\in C^\infty(\overline{B_1^+})$ and denote by  $\widetilde w$ its
fractional Kelvin transform defined as $\widetilde w(z)=|z|^{-(N-2s)}w\big(
\frac{z}{|z|^2}\big)$. Some computations (see also \cite{FW}) show that
\begin{align}\label{eq:12bis}
  &\int_{B_{1}^+}t^{1-2s}|\nabla w|^2dz+(N-2s)\int_{S_1^+}t^{1-2s}w^2dS
  =\int_{\R^{N+1}_+\setminus B_{1}^+}t^{1-2s}|\nabla \widetilde w|^2dz,\\
 \label{eq:10} &\int_{B_1'}|w|^{2^*(s)}dx=\int_{\R^{N}\setminus
    B_{1}'}|\widetilde w|^{2^*}dx.
\end{align}
In particular the function
$$
v(z)=
\begin{cases}
  w(z),&\text{if }z\in B_{1}^+,\\
\widetilde w(z) ,&\text{if }z\in \R^{N+1}_+\setminus  B_{1}^+
\end{cases}
$$ belongs to $\mathcal D^{1,2}(\R^{N+1}_+;t^{1-2s})$, so that, by
\eqref{eq:half_space_Sobolev} we have
\begin{equation}\label{eq:11}
    \bigg(\int_{\R^N}|v |^{2^*(s)}\,dx\bigg)^{\!\!\frac2{2^*(s)}}
\leq
\frac1{S_{N,s}\k_s}\int_{\R^{N+1}_+ }t^{1-2s}|\nabla v|^2dz.
\end{equation}
From \eqref{eq:12bis}, \eqref{eq:10}, and \eqref{eq:11}, it follows that
\begin{align*}
    \bigg(2 \int_{B_1'}|w|^{2^*(s)}dx\bigg)^{\!\!\frac2{2^*(s)}}
\leq
\frac2{S_{N,s}\k_s}\bigg[\int_{B_1^+ }t^{1-2s}|\nabla w|^2dz+\frac{N-2s}2\int_{S_1^+}t^{1-2s}w^2dS
\bigg]
\end{align*}
which, by density, yields the conclusion.
\end{pf}

\noindent Combining  Lemma  \ref{l:hardy_boundary2} and Lemma
\ref{l:sobolev_boundary} the following corollary follows.
\begin{Corollary}\label{c:1sobolev_boundary}
  For all $r>0$ and $w\in
  H^1(B_r^+;t^{1-2s})$, the following
  inequalities hold
  \begin{multline}\label{eq:14}
    \int_{B_r^+}t^{1-2s}|\nabla w|^2dz -\kappa_s
    \lambda\int_{B_r'}\frac{w^2}{|x|^{2s}}\,dx
    +\frac{N-2s}{2r}\int_{S_r^+}t^{1-2s}w^2dS\\
    \geq \kappa_s(\Lambda_{N,s}-\lambda) \int_{B_r'}\frac{w^2}{|x|^{2s}}\,dx
\end{multline}
and
  \begin{multline}\label{eq:15}
    \int_{B_r^+}t^{1-2s}|\nabla w|^2dz -\kappa_s
    \lambda\int_{B_r'}\frac{w^2}{|x|^{2s}}\,dx
    +\frac{N-2s}{2r}\int_{S_r^+}t^{1-2s}w^2dS\\
    \geq \frac{\Lambda_{N,s}-\lambda}{(1+\Lambda_{N,s})\widetilde S_{N,s}}
 \bigg(\int_{B_r'}| w |^{2^*(s)}\,dx\bigg)^{\!\!\frac2{2^*(s)}}.
\end{multline}
\end{Corollary}

\section{The Almgren type frequency
  function}\label{sec:monot-prop}
Let $R>0$ be such that $ B_R' \subset\subset\O$
 and $w\in  H^1(B_{R}^+;t^{1-2s})$ be a nontrivial solution to
 \begin{equation}\label{eq:wHextended}
\begin{cases}
    \dive(t^{1-2s}\nabla  w)=0,&\text{in } B_{R}^+,\\
-\lim_{t\to 0^+}t^{1-2s}\frac{\partial w}{\partial
  t}(t,x)=\kappa_s
\Big(\frac{\lambda}{|x|^{2s}}w+hw+f(x,w)\Big), &\text{on }B_R',
\end{cases}
\end{equation}
in a weak sense, i.e., for all $\varphi\in C^ \infty_c(B_{R}^+\cup
B_{R}')$,  we have that
\begin{equation}\label{eq:wH8}
 {\displaystyle{\int_{\R^{N+1}_+}t^{1-2s}\nabla w\cdot\nabla \varphi\,dt\,dx=
\kappa_s
\int_{B_R'}\bigg(
\dfrac{\lambda}{|x|^{2s}}w+hw+f(x,w)\bigg) \varphi\,dx}},
\end{equation}
   with $s,\lambda,h,f$ as in assumptions
\eqref{eq:ipo1},
\eqref{eq:ipoh}, and \eqref{eq:ipof}.

 For every $r\in (0,R]$ we define
\begin{equation}\label{D(r)}
  D(r)=\frac{1}{r^{N-2s}} \bigg[\int_{B_r^+}
t^{1-2s}|\nabla w|^2\,dt\,dx-
\kappa_s
\int_{B_r'}\bigg(
\dfrac{\lambda}{|x|^{2s}}w^2+hw^2+f(x,w)\bigg)\,dx\bigg]
\end{equation}
and
\begin{equation} \label{H(r)}
H(r)=\frac{1}{r^{N+1-2s}}\int_{S_r^+}t^{1-2s}w^2 \, dS
=
\int_{{\mathbb
      S}^{N}_+}\theta_1^{1-2s}w^2(r\theta)\,dS.
\end{equation}

The main result of this section is  the existence of  the limit as
$r\to 0^+$ of the \emph{Almgren's frequency} function (see \cite{GL}
and \cite{almgren}) associated to $w$
\begin{equation} \label{eq:31}
  {\mathcal N}(r)=\frac{D(r)}{H(r)}=\frac{\displaystyle r \left[\int_{B_r^+}
t^{1-2s}|\nabla w|^2\,dt\,dx-
\kappa_s
\int_{B_r'}\left(
\frac{\lambda}{|x|^{2s}}w^2+h(x)w^2+f(x,w)w\right)\,dx\right]}{\displaystyle\int_{S_r^+}t^{1-2s}w^2 \, dS}.
\end{equation}
We notice that, by Lemma \ref{l:hardy_boundary2}, $w(0,\cdot)\in L^2(B_{R}';|x|^{-2s})$ and so
 the
$L^1(0,R)$-function
$$
r\mapsto
\int_{S_r^+}
t^{1-2s}|\nabla w|^2\,dS,\quad\text{respectively}\quad
r\mapsto
\int_{\partial B_r'}\frac{w^2}{|x|^{2s}}\,{dS'},
$$
  is the weak derivative of the
$W^{1,1}(0,R)$-function
$$
r\to
\int_{B_r^+}
t^{1-2s}|\nabla w|^2\,dz,\quad\text{respectively}\quad
r\mapsto
\int_{B_r'}\frac{w^2}{|x|^{2s}}\,dx.
$$
In particular, for a.e. $r\in (0,R)$, $\frac{\partial
  w}{\partial\nu}\in L^2(S_r^+;t^{1-2s})$, where
  $\nu=\nu(z)=\frac{z}{|z|}$.\\
Next we observe the following integration by parts.
\begin{Lemma}\label{l:int_by_parts}
  For a.e. $r\in (0,R)$ and every $\widetilde\varphi\in
  C^\infty(\overline{B_r^+})$
$$
\int_{B_r^+}t^{1-2s}\nabla w\cdot\nabla\widetilde\varphi\,dz=
\int_{S_r^+}t^{1-2s}\frac{\partial w}{\partial\nu}\widetilde\varphi\,dS+
\kappa_s
\int_{B'_r}\bigg(
\frac{\lambda}{|x|^{2s}}w+hw+f(x,w)\bigg)\widetilde\varphi\,dx.
$$
\end{Lemma}

\begin{pf}
  It follows by testing \eqref{eq:wH8} with
  $\widetilde\varphi(z)\eta_n(|z|)$ where $\eta_n(\rho)=1$ if
  $\rho<r-\frac1n$, $\eta_n(r)=0$ if $\rho>r$, $\eta(\rho)=n(r-\rho)$
  if
  $r-\frac1n\leq \rho\leq r$, passing to the limit, and noticing that
  a.e. $r\in (0,R)$ is a Lebesgue point for the $L^1(0,R)$- function
  $r\mapsto \int_{S_r^+} t^{1-2s}\frac{\partial
    w}{\partial\nu}\widetilde\varphi\,dS$.
\end{pf}

\subsection{Regularity estimates and  Pohozaev-type identity }
In this section, we will prove local regularity estimates for a
general class of fractional elliptic equations in {Lemma}
\ref{lem:reg} below. This estimate will be useful for the blow-up
analysis and also for establishing  the Pohozaev identity in Theorem
\ref{t:pohozaev} below which is crucial in this paper. The proof of
Lemma \ref{lem:reg} uses mainly a result of Jin, Li and Xiong in
\cite{JLX} that we state here for sake of completeness.
\begin{Proposition}[\cite{JLX} Proposition 2.4]\label{prop:jlx} Let $v\in H^1(B_1^+; t^{1-2s})$ be a weak
solution to
$$
\begin{cases}
 \dive(t^{1-2s}\n v)=0,& \textrm{ in } B_{1}^+\\
-\lim_{t\to 0^+}t^{1-2s}v_t=c(x) v+b(x), &\textrm{on } B_{1}' ,
\end{cases}
$$
with $c,b\in L^p(B_1')$ for some $p>\frac{N}{2s}$. Then $v\in
L^\infty_{\rm loc}\left( B_1^+\cup B_1'\right)$ and  there exists
$C>0$ depending only on $N,s,p,\|c\|_{L^p(B_1')}$ such that
$$
\|v\|_{ L^\infty\left([0,1/2)\times B_{1/2}'\right) } \leq C \left(
\| v\|_{H^1(B_1^+; t^{1-2s})} + \|b\|_{L^p(B'_{1})} \right).
$$
Also $v\in C^{0,\a}_{\rm loc}\left( B_r^+\cup B_r'\right))$ for some
 $\a\in
(0,1)$ depending only on $N,s,p,\|c\|_{L^p(B_1')}$. In addition we
have
$$
\|v\|_{ C^{0,\a}\left([0,1/2)\times B_{1/2}'\right) } \leq C \left(
\| v\|_{L^\infty\left([0,1)\times B_{1}'\right)} +
\|b\|_{L^p(B'_{1})} \right).
$$
\end{Proposition}

We now state the following technical but crucial result.

\begin{Lemma}\label{lem:reg}
\begin{enumerate}[\rm (i)]
\item Let $r>0$ and $V\in L^q(B'_r)$ for some $q>\frac{N}{2s}$.
For every $t_0,r_0>0$ such that $[0,t_0) \times
  B_{r_0}'\Subset B_r^+\cup B_r'$ there exist  positive constants $A_1>0$, $\alpha\in(0,1)$
  depending  on $t_0,r_0,r,\|V\|_{L^q(B'_r)}$ such that for every
  $v\in H^1(B_{r}^+;t^{1-2s})$ solving
 \begin{equation*}
\begin{cases}
 \dive(t^{1-2s}\n v)=0, &\textrm{in } B_r^+\\
-\lim_{t\to 0^+}t^{1-2s}v_t=V(x)v,  &\textrm{on }B_r',
\end{cases}
\end{equation*}
we have that $v\in  C^{0,\a}([0,t_0) \times B_{r_0}')$ and
\be\label{eq:wCzaL}
\|v\|_{ C^{0,\a}([0,t_0) \times B_{r_0}') }\leq
 A_1 \|v\|_{H^1(B_{r}^+;t^{1-2s})}.
\ee
\item
Let $v\in H^1(B_{r}^+;t^{1-2s})\cap  C^{0,\a}(B_{r}^+)$ for some
$\alpha\in(0,1)$ be a function  satisfying
 \begin{equation*}
\label{eq:eqreg}
\begin{cases}
 \dive(t^{1-2s}\n v)=0, &\textrm{in } B_r^+\\
-\lim_{t\to 0^+}t^{1-2s}v_t=g(x,v),  &\textrm{on }B_r'  ,
\end{cases}
\end{equation*}
where $g\in C^1(B'_r\times\R)$ and
$|g(x,\rho)|\leq c(|\rho|+|\rho|^{p-1})$ for some $2<p\leq
2^*\!(s)=\frac{2N}{N-2s}$,  $c>0$, and every
$x\in B_r'$ and $\rho\in\R$. Let
 $t_0,r_0>0$ such that $[0,t_0) \times
  B_{r_0}'\Subset B_r^+$. Then there  exist   positive constants
$A_2,\b$ depending only on
$N$, $p$, $s$, $c$, $r$, $r_0$, $t_0$, $\| v\|_{H^1(B_{r}^+;t^{1-2s})
}$ and $ \|g\|_{C^{1}(B_{r_0}' \times[0,A_3])}$ where
$A_3=\|v\|_{ C^{0,\a}([0,t_0) \times B_{r_0}') }$,
 with $\b\in(0,1)$, such that
\be\label{eq:estnxwL}
\|\n_x v\|_{C^{0,\b}([0,t_0)\times B_{r_0}')}\leq  A_2 ,
\ee
\be\label{eq:estnxdtwL}
\|t^{1-2s} v_t\|_{C^{0,\b}([0,t_0)\times B_{r_0}')}\leq  A_2 .
\ee
\end{enumerate}
\end{Lemma}
\begin{remark}\label{rem:unif_reg}
  The dependence of the constant $A_2$ in Lemma \ref{lem:reg} on
  $\|v\|_{H^1(B_{r}^+;t^{1-2s}) }$ is continuous; in particular we can
  take the same $A_2$ for a family of solutions which are uniformly bounded in $H^1(B_{r}^+;t^{1-2s}) \cap  C^{0,\a}(B_{r}^+)$.
\end{remark}

\begin{pf}
Part (i) and \eqref{eq:wCzaL} follows from Proposition
\ref{prop:jlx}.

To prove (ii),  for $h\in\R^N$ with $|h|<<1$,
we set $v^h(t,x)=\frac{v(t,x+h)-v(t,x)}{|h|}$ for every
$x\in B_{r_0/2}'$. Then we have
$$
\begin{cases}
 \dive(t^{1-2s}\n v^h)=0,& \textrm{ in } B_{r_0/2}^+\\
-\lim_{t\to 0^+}t^{1-2s}v^h_t=c_h(x) v^h+b_h, &\textrm{on } B_{r_0/2}'  ,
\end{cases}
$$
where
$$
c_h(x)=\frac{g(x,v(0,x+h))-g(x,v(0,x))}{ v(0,x+h)-v(0,x)}\chi_{ \{v(0,x+h)\neq v(0,x) \}}(x)
$$
with $\chi_A$ being the characteristic function of a set $A$ and
$$
b_h(x)=\frac{g(x+h,v(0,x+h))-g(x,v(0,x+h))}{|h|}.
$$
Let $A_3=\|v\|_{ C^{0,\a}([0,t_0) \times B_{r_0}') }$. Then we have
$$
\|c_h\|_{L^\infty(B_{r_0/4}')} + \|b_h\|_{L^\infty(B_{r_0/4}')} \leq
 \| g_\rho\|_{L^\infty(B_{r_0/2}'\times[0,A_3])}+
\| \n_x g\|_{L^\infty(B_{r_0/2}'\times[0,A_3]),
}$$
for every small $h$.
 Applying once again  Proposition
\ref{prop:jlx},
\begin{eqnarray*}
\|v^h\|_{ C^{0,\a}([0,t_0/8)\times B_{r_0/8}') }
&\leq& C ( \| v^h\|_{L^\infty([0,t_0/4)\times B_{r_0/4}')} + \|b_h\|_{L^\infty(B'_{r_0/2})} )
\\
&\leq &C  \|v^h\|_{L^2([0,t_0/2)\times B_{r_0/2}' ;t^{1-2s })} +C \|\n g\|_{L^\infty(B_{r_0/2}'\times[0,A_3])} \\
&\leq &C  \|\n v\|_{L^2 ([0,t_0) \times B_{r_0}';t^{1-2s })} +C \|\n g\|_{L^\infty(B_{r_0/2}'\times[0,A_3])}\\
&\leq& A_2
\end{eqnarray*}
for every small $h$.   By the Arzelà-Ascoli Theorem, passing to the
limit as $|h|\to 0$, we conclude that $z\mapsto \n_x v(z)\in C([0,t_0/8)
\times B_{r_0/8}')$ and estimate \eqref{eq:estnxwL} holds
for all $\b\in(0,\a)$.

Since the map $x\mapsto g(x,v(0,x))\in C^{0,\b}( B_{r_0/8}')$, estimate \eqref{eq:estnxdtwL}
follows from [\cite{CaSi}, Lemma 4.5].~\end{pf}

Because of the presence of the  Hardy potential $|x|^{-2s}$, we
cannot expect the solution $w$ of \eqref{eq:wHextended} to be
$L^q_{\rm loc}$ near the origin for every $q$, but we can expect $w$ to
be in a space better than $L^{2^*(s)}_{\rm loc}$. This is what we will
prove in the next result. 

\begin{Lemma}\label{lem-upCrit}
Let  $w$ be a solution to \eqref{eq:wHextended} in the sense of \eqref{eq:wH8}. Then
 there exist $p_0>2^*(s)$ and $R_0\in(0,R)$
such that $ w\in L^{p_0}(B_{R_0}^+)$.
\end{Lemma}
\begin{pf}
By \eqref{eq:ipoh}, there are $\d>0$ and  $R_\d\in(0,R)$ such that
\be\label{eq:hlesq}
(\l+|x|^{2s}|h(x)|)\leq\l+\d < \L_{N,s},\quad\text{for all }x\in B'_{R_\d}.
\ee
 Let $\b>1$.
For all $L>0$,  we define $F_L(\t)=|\t|^\b$ if $|\t|<L$ and $F_L(\t)=\b L^{\b-1} |\t|+(1-\b)L^\b$ if $|\t|\geq L$.
 Put $G_=\frac1\beta F_LF_L'$.
It is easy to verify that, for all $\tau\in\R$,
\be\label{eq:G_L}
\tau G_L(\t)\leq \t^2 G_L'(\t),\quad \t G_L(\t)\leq (F_L(\t))^2,\quad (F_L'(\t))^2\leq \b G_L'(\t).
 \ee
Let  $\eta \in C^\infty_c(B^+_{R_\delta}\cup B'_{R_\delta})$ be a
radial cut-off function
such that $\eta\equiv 1$ in $B^+_{R_\delta/2}$.
It is  clear that $\zeta:=\eta^2 G_L(w)\in H^1(B_R^+;t^{1-2s})$ and
$ F_L(w)\in H^1(B_R^+;t^ {1-2s})$. Using $\zeta$ as a test in
\eqref{eq:wHextended}, from \eqref{eq:hlesq} and integration by parts,
we have that
\begin{align*}
\int_{B_{R_\delta}^+}&t^{1-2s}\eta^2|\n w|^2G_L'(w)dtdx
-\k_s(\l+\d)\int_{B_{R_\delta}'}|x|^{-2s}\eta^ 2w G_L(w)dx\\
&\leq -2 \int_{B_{R_\delta}^+}t^{1-2s}\eta\n w\cdot\n\eta G_L(w)dtdx \\
&\quad+c \int_{B_{R_\delta} '}\eta^2 w G_L(w)dx+c \int_{B_{R_\delta} '}\eta^2|w|^{2^*(s)-2}w G_L(w)dx,
\end{align*}
for some positive $c>0$ depending only on $C_f,s,p,N$.
By Young's inequality and \eqref{eq:G_L}, we have that, for every $\sigma>0$,
$$
\left|2(\eta\n w)\cdot(w\n\eta)\frac{ G_L(w)}{w}\right| \leq
\frac{\sigma}{2} \eta^2|\n w|^2G_L'(w)+\frac{2}{\sigma}|\n\eta|^2(F_L(w))^2,
$$
hence we obtain that
\begin{multline*}
\left(1-\frac{\sigma}{2}\right)\int_{B_{R_\delta}^+}t^{1-2s}\eta^2|\n w|^2G_L'(w)\,dt\,dx
-\k_s(\l+\d)\int_{B_{R_\delta}'}|x|^{-2s}\eta^ 2w G_L(w)\,dx\\
\leq  c \int_{B_{R_\delta}'}\eta^2(|w|^{2^*(s)-2}+1)w\, G_L(w)\,dx
+\frac2\sigma \int_{ B_{R_\delta}^+}t^{1-2s} |\n \eta|^2\psi^2 dt\,dx,
\end{multline*}
where we have set $\psi=F_L(w)$.
Therefore by \eqref{eq:G_L}
\begin{multline*}
\frac{1}{\b}\left(1-\frac{\sigma}{2}\right)\int_{B_{R_\delta}^+}t^{1-2s}\eta^2|\n \psi |^2dt\,dx
-\k_s(\l+\d) \int_{ B_{R_\delta}'}|x|^{-2s}(\eta\psi)^2dx\\
\leq
c  \int_{B_{R_\delta}'}(|w|^{2^*(s)-2}+1)(\eta\psi)^2\,dx
+\frac2\sigma \int_{ B_{R_\delta}^+}t^{1-2s} |\n \eta|^2\psi^2 dt\,dx.
\end{multline*}
Since $ |\n(\eta \psi) |^2\leq (1+\sigma)\eta^2|\n \psi|^2+(1+\frac1\sigma)
\psi^2|\n \eta|^2$, we have that
\begin{multline*}
\frac{1}{\b(1+\sigma)}\left(1-\frac{\sigma}{2}\right)\int_{B_{R_\delta}^+}t^{1-2s}|\n
(\eta \psi) |^2dt\,dx
-\k_s(\l+\d) \int_{ B_{R_\delta}'}|x|^{-2s}(\eta\psi)^2dx\\
\leq
 c \int_{B_{R_\delta}'}(|w|^{2^*(s)-2}+1)(\eta\psi)^2\,dx
+ C(c,\beta,\sigma,R_\delta)
\int_{ B_{R_\delta}^+}t^{1-2s} \psi^2 dt\,dx
\end{multline*}
for some positive $C(c,\beta,\sigma,R_\delta)>0$ depending only on $c$, $\sigma$,
$R_\delta$, and $\beta$.
By H\"{o}lder inequality and Lemma \ref{l:sobolev_boundary}, we have  that
\begin{multline*}
  \int_{B_{R_\delta}'}(|w|^{2^*(s)-2}+1)(\eta\psi)^2\,dx\\
  \leq \widetilde S_{N,s}
  \left[\bigg(\int_{B_{R_\delta}'}|w|^{2^*(s)}dx\bigg)
    ^{\!\!\frac{{2^*(s)}-2}{{2^*(s)}}}+|B_{R_\delta}'|^{\frac{2s}{N}}
  \right] \int_{B_{R_\delta}^+}t^{1-2s}|\n (\eta\psi) |^2dt\,dx.
\end{multline*}
We deduce that
\be\label{eq:mtial}
A\,\int_{B_{R_\delta}^+}t^{1-2s}|\n (\eta\psi) |^2dt\,dx-\k_s(\l+\d) \int_{ B_{R_\delta}'}|x|^{-2s}(\eta\psi)^2dx\leq
{\rm const\,} \int_{ B_{R_\delta}^+}t^{1-2s} \psi^2 dtdx,
\ee
 for some positive ${\rm const}>0$ depending only on $C_f,s,p,N, \beta,\sigma,R_\delta$,
where
$$
A=\frac{1}{(1+\sigma)\b}\left(1-\frac{\sigma}{2}\right)- c \widetilde S_{N,s}
  \left[\bigg(\int_{B_{R_\delta}'}|w|^{2^*(s)}dx\bigg)
    ^{\!\!\frac{{2^*(s)}-2}{{2^*(s)}}}+|B_{R_\delta}'|^{\frac{2s}{N}}
  \right].
$$
From Hardy inequality \eqref{eq:half_space_hardy}, we have that
\begin{multline*}
A\,\int_{B_{R_\delta}^+}t^{1-2s}|\n (\eta\psi) |^2dt\,dx-\k_s(\l+\d) \int_{ B_{R_\delta}'}|x|^{-2s}(\eta\psi)^2dx\\\geq
\bigg(A-\frac{\lambda+\delta}{\Lambda_{N,s}}\bigg)
\int_{B_{R_\delta}^+}t^{1-2s}|\n (\eta\psi) |^2dt\,dx,
\end{multline*}
and, by  \eqref{eq:hlesq}, we can choose $\b$ sufficiently close to
$1$ and $\sigma,R_\delta$ sufficiently small such  that
$$
A-\frac{\lambda+\delta}{\Lambda_{N,s}}>0.
$$
Hence we have that
\begin{equation*}
C\int_{B_{R_\delta}^+}t^{1-2s}|\n (\eta\psi) |^2dt\,dx\leq \int_{ B_{R_\delta}^+}t^{1-2s} \psi^2 dtdx,
\end{equation*}
for some constant $C>0$ depending on $f,h,s,p,N,\beta,\e,w,\lambda,\delta$.
From Lemma \ref{l:sobolev_boundary} it follows that
$$
C\widetilde S_{N,s}^{-1}\bigg(\int_{B_{R_\delta}'}| \eta F_L(w) |^{2^*(s)}dx\bigg)^{\!\!\frac{2}{2^*(s)}}\leq
  \int_{B_{R_\delta}^+}t^{1-2s}|w|^{2\b}dtdx\quad\text{for all }L\geq0.
$$
Hence by taking the limit as  $L \to +\infty$
\be\label{eq:eFLw}
C\widetilde S_{N,s}^{-1}\bigg(\int_{B_{R_\delta/2}'}|w|^{\beta 2^*(s)}dx\bigg)^{\!\!\frac{2}{2^*(s)}}\leq
 \int_{B_{R_\delta}^+}t^{1-2s}|w|^{2\b}dtdx.
\ee
The conclusion follows since $\b>1$ and $H^1(B_R^+;t^{1-2s})\hookrightarrow L^q(B_R^+;t^{1-2s})$
for some $q>2$, see for instance [\cite{FKS}, Theorem 1.2].
\end{pf}

\begin{remark}\label{rem:rego}
From Lemma \ref{lem-upCrit}, we deduce that, if
$w\in  H^1(B_{R}^+;t^{1-2s})$ is a weak  solution to
\eqref{eq:wHextended}, then
$\frac{\lambda}{|x|^{2s}}+h+\frac{f(x,w)}w\in L^q_{\rm
  loc}(B'_{R}\setminus\{0\})$ for some $\frac{N}{2s}<q\leq \frac{p_0}{p-2}$
and hence, from
Lemma \ref{lem:reg}, we conclude that $w\in  C^{0,\a}_{\rm
  loc}(\overline{B_{r}^+}\setminus\{0\})$,
$\n_x v\in C^{0,\b}_{\rm
  loc}(\overline{B_{r}^+}\setminus\{0\})$, and $t^{1-2s} v_t\in C^{0,\b}_{\rm
  loc}(\overline{B_{r}^+}\setminus\{0\})$ for all $r\in(0,R)$ and some $0<\beta<\alpha<1$.
\end{remark}

\noindent We next prove the following Pohozaev-type identity which
will be used to compute $D'$ (see Lemma \ref{l:dprime} below) and
therefore $\calN'$.

\begin{Theorem} \label{t:pohozaev}
Let $w$ solves \eqref{eq:wH8}. Then
for a.e. $r\in (0,R)$ there holds
\begin{multline}\label{eq:poho}
  -\frac{N-2s}2\bigg[\int_{B_r^+} t^{1-2s}|\nabla w|^2dz
  -\kappa_s\lambda
  \int_{B_r'}\frac{w^2}{|x|^{2s}}dx\bigg]\\+\frac{r}{2}\bigg[
\int_{S_r^+} t^{1-2s}|\nabla w|^2dS -\kappa_s\lambda
  \int_{\partial B_r'}\frac{w^2}{|x|^{2s}}dS'\bigg]\\
  =r\int_{S_r^+}t^{1-2s}\bigg|\frac{\partial w}{\partial
    \nu}\bigg|^2\,dS
  -\frac{\kappa_s}2\int_{B_r'}(Nh+\nabla h\cdot x)w^2\,dx
 +
  \frac{r\kappa_s}2\int_{\partial B_r'}hw^2\,dS' \\
+r\kappa_s\int_{\partial B_r'} F(x,w)\, dS'-\kappa_s\int_{B_r'}
[\nabla_xF(x,w)\cdot x+NF(x,w)]\, dx
\end{multline}
and
\begin{multline}\label{eq:poho2}
 \int_{B_r^+} t^{1-2s}|\nabla w|^2dz
  -\kappa_s\lambda
  \int_{B_r'}\frac{w^2}{|x|^{2s}}dx\\=
\int_{S_r^+}t^{1-2s}\frac{\partial w}{\partial\nu}w\,dS+
\kappa_s
\int_{B'_r}\bigg(hw^2+f(x,w)w\bigg)\,dx.
  \end{multline}
\end{Theorem}

\begin{pf}
We write our problem in the form
$$\begin{cases}
 \dive(t^{1-2s}\n w)=0, & \textrm{ in } B_R^+,\\
-\lim_{t\to 0^+}t^{1-2s}w_t=G(x,w), & \textrm{on } B_R',
\end{cases}
$$
where $G\in C^1(B_R'\setminus\{0\}\times \R)$,
$G(x,\varrho)=\kappa_s
\big(\frac{\lambda}{|x|^{2s}}\varrho+h(x)\varrho+f(x,\varrho)\big)$.

We have, on $B_R^+$, the formula
 \be\label{eq:diver} \dive\left(
  \frac{1}{2}t^{1-2s}|\n w|^2 z-t^{1-2s}(z\cdot \n w) \n w
\right)=\frac{N-2s}{2} t^{1-2s}|\n w|^2-(z\cdot \n w )
\dive(t^{1-2s}\n w).  \ee
Let $\rho<r<R$. Now we integrate by parts
over the set $O_\delta:= (B_r^+\setminus{\ov{B_\rho^+}})\cap\{(t,x),\,
t>\delta\} $ with $\delta>0$. We have
\begin{align*}
\frac{N-2s}{2}& \int_{O_\delta} t^{1-2s}|\nabla w(z)|^2dz =
-\frac{1}{2}\delta^{2-2s}\int_{B_{\sqrt{r^2-\delta^2}}' \setminus
  B_{\sqrt{\rho^2-\delta^2}}'}|\n w|^2(\delta,x )dx\\
&+ \delta^{2-2s}\int_{B_{\sqrt{r^2-\delta^2}}'  \setminus
  B_{\sqrt{\rho^2-\delta^2}}' }| w_t|^2(\delta,x )dx \hspace{2cm}\\
&+\frac{r}{2}\int_{S_r^+\cap\{t>\delta\}} t^{1-2s}|\nabla w|^2dS-
r\int_{S_r^+\cap\{t>\delta\}}t^{1-2s}\bigg|\frac{\partial w}{\partial\nu}\bigg|^2dS\\
&-\frac{\rho}{2}\int_{S_\rho^+\cap\{t>\delta\}} t^{1-2s}|\nabla w|^2dS
+\rho\int_{S_\rho^+\cap\{t>\delta\}}t^{1-2s}\bigg|\frac{\partial w}{\partial\nu}\bigg|^2dS\\
&+ \int_{B_{\sqrt{r^2-\delta^2}}'  \setminus
  B_{\sqrt{\rho^2-\delta^2}}'}(x\cdot \n_x w(\delta,x)) \,\delta^{1-2s}w_t(\delta,x)\,dx.
\end{align*}
We now claim that there exists a sequence $\delta_n\to 0 $ such that
$$
\lim_{n\to \infty} \left[  \frac{1}{2}\delta_n^{2-2s}\int_{B_r' }|\n w|^2(\delta_n,x )dx+
 \delta_n^{2-2s}\int_{B_r' }| w_t|^2(\delta_n,x )dx \right]=0.
$$
If no such sequence exists, we  would have
$$
\liminf_{\delta\to 0 } \left[  \frac{1}{2}\delta^{2-2s}
\int_{B_r' }|\n w|^2(\delta,x )dx+ \delta^{2-2s}\int_{B_r' }| w_t|^2(\delta,x )dx \right]\geq C>0
$$
and thus there exists $\delta_0>0$ such that
$$
 \frac{1}{2}\delta^{2-2s}\int_{B_r' }|\n w|^2(\delta,x )dx+ \delta^{2-2s}\int_{B_r' }| w_t|^2(\delta,x )dx
\geq \frac{C}{2}\quad \text{for all } \delta\in(0,\delta_0).
$$
It follows that
$$
 \frac{1}{2}\delta^{1-2s}\int_{B_r' }|\n w|^2(\delta,x )dx+ \delta^{1-2s}\int_{B_r' }| w_t|^2(\delta,x )dx
 \geq \frac{C}{2\delta}\quad \text{for all } \delta\in(0,\delta_0)
$$
and so integrating the above inequality on $(0,\delta_0)$ we contradict
the fact that $w\in H^1(B_{R}^+;t^{1-2s})$.

Next, from the Dominated Convergence Theorem, Lemma \ref{lem:reg} and
Remark \ref{rem:rego}, we have that
$$
\lim_{\delta\to 0} \int_{B_{\sqrt{r^2-\delta^2}}'  \setminus
  B_{\sqrt{\rho^2-\delta^2}}'}(x\cdot \n_x w(\delta,x))
\,\delta^{1-2s}w_t(\delta,x)\,dx= -\int_{B_r'\setminus B_\rho' }(x\cdot \n_x w)
\,G(x,w)\,dx.
$$
We conclude that (replacing $O_\delta$ with $O_{\delta_n}$, for a sequence $\delta_n\to 0$) that
\begin{multline}\label{eq:pohozG}
\frac{N-2s}2 \int_{B_r^+\setminus{{B_\rho^+}}} t^{1-2s}|\nabla w(z)|^2dz =
\frac{r}{2}\int_{S_r^+} t^{1-2s}|\nabla w|^2dS-r\int_{S_r^+}t^{1-2s}\bigg|\frac{\partial w}{\partial
    \nu}\bigg|^2\,dS\\
-\frac{\rho}{2}\int_{S_\rho^+} t^{1-2s}|\nabla w|^2dS-\rho\int_{S_\rho^+}t^{1-2s}\bigg|\frac{\partial w}{\partial
    \nu}\bigg|^2\,dS
- \int_{B_r'\setminus B_\rho' }(x\cdot \n_x w) \,G(x,w)\,dx.
\end{multline}
Furthermore, integration by parts yields
\begin{align}\label{eq:9}
  \int_{B_r'\setminus B_\rho' }&(x\cdot \n_x w)
\,G(x,w)\,dx=-\frac{N-2s}2\kappa_s\lambda
  \int_{B_r'\setminus B_\rho'}\frac{w^2}{|x|^{2s}}dx\\
\notag&-\frac{\kappa_s}2\int_{B_r'\setminus B_\rho'}(Nh(x)+\nabla h(x)\cdot x)w^2\,dx
+\kappa_s\lambda\frac{r}{2}
  \int_{\partial B_r'}\frac{w^2}{|x|^{2s}}{dS'}\\
\notag&+
  \frac{r\kappa_s}2\int_{\partial B_r'}h(x)w^2\,{dS'}-\kappa_s\lambda\frac{\rho}{2}
  \int_{\partial B_\rho'}\frac{w^2}{|x|^{2s}}{dS'}
-
  \frac{\rho\kappa_s}2\int_{\partial B_\rho'}h(x)w^2\,{dS'}\\
\notag&-\kappa_s\int_{ B_r'\setminus B_\rho'}
[\nabla_xF(x,w)\cdot x+NF(x,w)]\, dx\\
\notag&+r\kappa_s\int_{\partial B_r'} F(x,w)\, {dS'}
-\rho\kappa_s\int_{\partial B_\rho'} F(x,w)\, {dS'}.
\end{align}
 Since $w\in H^1(B_{R}^+;t^{1-2s})$, in view of Lemma
 \ref{l:hardy_boundary2} and \eqref{eq:sobolev}, there
 exists a sequence $\rho_n\to 0 $ such that
$$
\lim_{n\to \infty} \rho_n\bigg[\int_{S_{\rho_n}^+} t^{1-2s}|\nabla
w|^2dS + \int_{\partial
  B_{\rho_n}'}\frac{w^2}{|x|^{2s}}dS'+\int_{\partial B_{\rho_n}'}
|F(x,w)|\, {dS'} \bigg] =0.
$$
Hence, taking $\rho=\rho_n$ and letting $n\to\infty$ in
\eqref{eq:pohozG} and \eqref{eq:9}, we obtain \eqref{eq:poho}.

\eqref{eq:poho2} follows from Lemma \ref{l:int_by_parts} and density
of $C^\infty(\overline{B_r^+})$ in $H^1(B_r^+;t^{1-2s})$.
\end{pf}

\subsection{On the Almgren type frequency $\calN$}
In this section, we shall study the differentiability of $\calN$,
it's limit at 0 and provide estimates of $\calN'$.
\begin{Lemma} \label{l:hprime}
$H\in C^1(0,R)$ and
\begin{align}\label{H'}
 & H'(r)=\frac{2}{r^{N+1-2s}} \int_{S_r^+}t^{1-2s}w\frac{\partial
    w}{\partial\nu}
  \, dS, \quad \text{for every } r\in (0, R),\\
 &\label{H'2} H'(r)=\frac2r D(r), \quad \text{for every } r\in (0, R).
\end{align}
\end{Lemma}

\begin{pf}
Fix $r_0\in (0,R)$ and consider the limit
\begin{equation} \label{limite} \lim_{r\rightarrow r_0}
  \frac{H(r)-H(r_0)}{r-r_0} = \lim_{r\rightarrow r_0} \int_{{\mathbb
      S}^{N}_+}\theta_1^{1-2s} \frac{|w(r\theta)|^2-|w(r_0\theta)|^2}{r-r_0}dS.
\end{equation}
Since $w\in C^1(\R^{N+1}_+)$, then, for every $\theta\in {\mathbb
  S}^{N}_+$,
\begin{equation}\label{55} \lim_{r\rightarrow r_0}
  \frac{|w(r\theta)|^2-|w(r_0\theta)|^2}{r-r_0}
  =2\frac{\partial w}{\partial
      \nu}(r_0\theta)\,w(r_0\theta).
\end{equation}
On the other hand, for any $r\in (r_0/2,R)$ and $\theta\in {\mathbb
  S}^{N}_+$ we have
$$
\left| \frac{|w(r\theta)|^2-|w(r_0\theta)|^2}{r-r_0} \right| \leq 2
\sup_{B_{ R}^+\setminus B^+_{\frac{r_0}{2}}} |w| \cdot \sup_{B_{
    R}^+\setminus B^+_{\frac{r_0}{2}} } \bigg|\frac{\partial
  w}{\partial \nu}\bigg|\leq
2 \sup_{B_{ R}^+\setminus
  B^+_{\frac{r_0}{2}}} |w| \cdot \sup_{B_{ R}^+\setminus
  B^+_{\frac{r_0}{2}} } \bigg(\bigg|
\frac{t}{|z|}w_t\bigg|+\bigg|\frac{\nabla_xw\cdot x}{|z|}\bigg|\bigg)$$
and hence, by (\ref{limite}), (\ref{55}), Lemma \ref{lem:reg} and the Dominated
Convergence Theorem, we obtain that
$$
H'(r_0)=
2\int_{{\mathbb
      S}^{N}_+}\theta_1^{1-2s} \frac{\partial w}{\partial
      \nu}(r_0\theta)\,w(r_0\theta)dS(\theta)
=
\frac{2}{r_0^{N+1-2s}} \int_{S_{r_0}^+}t^{1-2s}w\frac{\partial
    w}{\partial\nu}
  \, dS.
$$
The continuity of $H'$ on the interval $(0,R)$ follows by the
representation of $H'$ given above, Lemma \ref{lem:reg}, and the
Dominated Convergence Theorem.

Finally, \eqref{H'2} follows from \eqref{H'}, \eqref{D(r)}, and \eqref{eq:poho2}.
\end{pf}

 The regularity of the
function $D$ is established in the following lemma.

\begin{Lemma}\label{l:dprime}
  The function $D$ defined in (\ref{D(r)}) belongs to $W^{1,1}_{{\rm\
      loc}}(0, R)$ and
\begin{align}\label{D'F}
  D'(r)=&\frac{2}{r^{N+1-2s}} \bigg[ r\int_{S_r^+}
  t^{1-2s}\left|\frac{\partial w}{\partial \nu}\right|^2 dS
  -\kappa_s\int_{B_r'}\Big(sh+\frac12(\nabla h\cdot x)\Big)w^2\,dx\bigg]\\
  \notag & \ +\frac{\kappa_s}{r^{N+1-2s}}\int_{B_r'}
  \big((N-2s)f(x,w)w-2NF(x,w)-
  2\nabla_xF(x,w)\cdot x\big)\, dx \\
  \notag & \ +\frac{\kappa_s}{r^{N-2s}}\int_{\partial B_r'}
  \big(2F(x,w)-f(x,w)w\big)\, {dS'}
\end{align}
in a distributional sense and for a.e. $r\in (0,R)$.
\end{Lemma}

\begin{pf}
For any $r\in (0,r_0)$ let
\begin{align}\label{I(r)}
I(r)= \int_{B_r^+}
t^{1-2s}|\nabla w|^2\,dt\,dx-
\kappa_s
\int_{B_r'}\bigg(
\dfrac{\lambda}{|x|^{2s}}w^2+hw^2+f(x,w)w\bigg)\,dx.
\end{align}
From
the fact that $w\in  H^1(B_R^+;t^{1-2s})$, Lemma
 \ref{l:hardy_boundary2}, and \eqref{eq:sobolev},  we deduce that $I\in
 W^{1,1}(0,R)$ and
\begin{equation} \label{I'(r)}
I'(r) =  \int_{S_r^+}
t^{1-2s}|\nabla w|^2\,dS-
\kappa_s
\int_{\partial B_r'}\bigg(
\dfrac{\lambda}{|x|^{2s}}w^2+hw^2+f(x,w)w\bigg)\,{dS'}
\end{equation}
 for a.e. $r\in (0,R)$ and in the distributional sense.
Therefore $D\in W^{1,1}_{{\rm loc}}(0,R)$ and, using
\eqref{eq:poho}, (\ref{I(r)}), and (\ref{I'(r)}) into
\begin{align*}
  D'(r)=r^{2s-1-N}[-(N-2s)I(r)+rI'(r)],
\end{align*}
we obtain (\ref{D'F}) for a.e. $r\in (0,R)$ and in the distributional sense.
\end{pf}

\noindent
Before going on, we recall that $w$ is nontrivial and satisfies
\eqref{eq:wH8}.
We prove now that, if $w\not\equiv 0$,
 $H(r)$ does not vanish for $r$ sufficiently small.

 \begin{Lemma} \label{welld}
There exists $R_0\in(0,R)$   such that
$H(r)>0$ for any $r\in (0,R_0)$, where $H$ is defined by (\ref{H(r)}).
\end{Lemma}

\begin{pf}
Clearly from assumption \eqref{eq:ipo1}, there exists $R_0\in (0,R)$ such that
\begin{align}\label{eq:r_0}
  \frac{\lambda}{\Lambda_{N,s}}+\frac{C_hR_0^\e}{\Lambda_{N,s}}
+C_f S_{N,s}^{-1}\big(\tfrac{\omega_{N-1}}N\big)^{\!\frac{2^*(s)-p}{2^*(s)}}R_0^{\frac{N(2^*(s)-p)}{2^*(s)}}
\|w\|_{L^{2^*(s)}(B_{R_0}')}^{p-2}<1 ,
\end{align}
where $\omega_{N-1}$ denotes the volume of the unit sphere ${\mathbb
  S}^{N-1}$, i.e. $\omega_{N-1}=\int_{{\mathbb S}^{N-1}}dS$.

Next suppose by contradiction that there exists $r_0\in(0,R_0)$ such that
$H(r_0)=0$. Then  $w= 0$ a.e. on $S_{r_0}^+$. From  \eqref{eq:poho2} it follows that
\begin{align*}
 \int_{B_{r_0}^+} t^{1-2s}|\nabla w|^2dz
  -\kappa_s\lambda
  \int_{B_{r_0}'}\frac{w^2}{|x|^{2s}}dx-
\kappa_s
\int_{B'_{r_0}}\bigg(
h(x)w^2+f(x,w)w\bigg)\,dx=0.
\end{align*}
From Lemma \ref{l:hardy_boundary2}, assumptions \eqref{eq:ipoh}-\eqref{eq:ipof},
H\"older's inequality, and \eqref{eq:sobolev}, it follows that
\begin{align*}
0&=\int_{B_{r_0}^+} t^{1-2s}|\nabla w|^2dz
  -\kappa_s\lambda
  \int_{B_{r_0}'}\frac{w^2}{|x|^{2s}}dx-
\kappa_s
\int_{B'_{r_0}}\bigg(
h(x)w^2+f(x,w)w\bigg)\,dx\\
\notag & \geq
  \bigg[1-
  \frac{\lambda}{\Lambda_{N,s}}-\frac{C_hR_0^\e}{\Lambda_{N,s}}-C_f S_{N,s}^{-1}
  \big(\tfrac{\omega_{N-1}}N\big)^{\!\frac{2^*(s)-p}{2^*(s)}}R_0^{\frac{N(2^*(s)-p)}{2^*(s)}}
\|w\|_{L^{2^*(s)}(B_{R_0}')}^{p-2} \bigg] \int_{B_{r_0}^+} t^{1-2s}|\nabla w|^2dz
,
\end{align*}
which, together with \eqref{eq:r_0}, implies $w\equiv 0$ in
$B_{r_0}^+$ by Lemma \ref{l:hardy_boundary}.  Classical
unique continuation principles for second order elliptic equations
with locally bounded coefficients (see e.g.  \cite{wolff}) allow to
conclude that $w=0$ a.e. in $B_R^+$, a contradiction.
\end{pf}

\noindent Letting $R_0$ be as in Lemma \ref{welld} and  recalling
\eqref{eq:31}, the \emph{Almgren type frequency}
  function
\begin{equation}\label{N(r)}
{\mathcal N}(r)=\frac{D(r)}{H(r)}
\end{equation}
is well defined in $(0,R_0)$.
 Using Lemmas \ref{l:hprime} and \ref{l:dprime}, we can now
 compute the derivative of ${\mathcal N}$.

 \begin{Lemma}\label{mono} The function
   ${\mathcal N}$ defined in (\ref{N(r)}) belongs to $W^{1,1}_{{\rm
       loc}}(0, R_0)$ and
\begin{align}\label{formulona}
{\mathcal N}'(r)=\nu_1(r)+\nu_2(r)
\end{align}
in a distributional sense and for a.e. $r\in (0,R_0)$,
where
\begin{align}\label{eq:nu1}
\nu_1(r)=&\frac{2r\Big[
    \left(\int_{S_r^+}
  t^{1-2s}\left|\frac{\partial w}{\partial \nu}\right|^2 dS\right) \cdot
    \left(\int_{S_r^+}
  t^{1-2s}w^2\,dS\right)-\left(
\int_{S_r^+}
  t^{1-2s}w\frac{\partial w}{\partial \nu}\, dS\right)^{\!2} \Big]}
{\left(
\int_{S_r^+}
  t^{1-2s}w^2\,dS\right)^2},
\end{align}
$\nu_1\geq0$, and
\begin{align}\label{eq:nu2}
  \nu_2(r)= & -\kappa_s\frac{\int_{B_r'} (2sh+\nabla h\cdot
    x)|w|^2\,dx}{\int_{S_r^+} t^{1-2s}w^2\,dS}+\kappa_s\frac{r\int_{\partial
      B_r'}\big(2F(x,w)-f(x,w)w\big)\, dS'} {\int_{S_r^+}
    t^{1-2s}w^2\,dS}
  \\[10pt]
  &\notag +\kappa_s\frac{\int_{B_r'}
    \big((N-2s)f(x,w)w-2NF(x,w) -2\nabla_xF(x,w)\cdot
    x\big)\, dx} {\int_{S_r^+} t^{1-2s}w^2\,dS}.
\end{align}
\end{Lemma}

\begin{pf} From Lemmas \ref{l:hprime}, \ref{welld}, and
  \ref{l:dprime}, it follows that ${\mathcal N}\in W^{1,1}_{{\rm
      loc}}(0,R_0)$. From
\eqref{H'2} it follows that
$$
{\mathcal N}'(r)=\frac{D'(r)H(r)-D(r)H'(r)}{(H(r))^2}
=\frac{D'(r)H(r)-\frac{1}{2} r (H'(r))^2}{(H(r))^2}
$$
and the proof of the lemma easily follows from (\ref{H'}) and
(\ref{D'F}). Now it is easy to see that $\nu_1\geq0$ by Schwarz's
inequality.
\end{pf}

\noindent We now prove that ${\mathcal N}(r)$ admits a finite limit as
$r\to 0^+$. To this aim, the following estimate plays a crucial role.

\begin{Lemma} \label{l:stimasotto}
 Let   ${\mathcal N}$ be the function defined in (\ref{N(r)}).
There
exist $\tilde R\in (0,R_0)$ and a  constant
$\overline{C}>0$  such that
\begin{multline}\label{eq:13}
\int_{B_r^+}
t^{1-2s}|\nabla w|^2\,dt\,dx-
\kappa_s
\int_{B_r'}\bigg(
\dfrac{\lambda}{|x|^{2s}}w^2+hw^2+f(x,w)w\bigg)\,dx\\\geq
-\bigg(\frac{N-2s}{2r}\bigg)\int_{S_r^+}t^{1-2s}w^2dS
+\overline{C}\left(\int_{B_r'}\dfrac{w^2}{|x|^{2s}}\,dx+\bigg(\int_{B_r'}|w|^{2^*(s)}\,dx
\bigg)^{\!\!\frac2{2^*(s)}}
\right),
\end{multline}
\begin{multline}\label{eq:13bis}
\int_{B_r^+}
t^{1-2s}|\nabla w|^2\,dt\,dx-
\kappa_s
\int_{B_r'}\bigg(
\dfrac{\lambda}{|x|^{2s}}w^2+hw^2+f(x,w)w\bigg)\,dx\\\geq
-\bigg(\frac{N-2s}{2r}\bigg)\int_{S_r^+}t^{1-2s}w^2dS
+\overline{C}\int_{B_r^+}
t^{1-2s}|\nabla w|^2\,dt\,dx,
\end{multline}
and
\begin{equation}\label{Nbelow}
   {\mathcal N}(r)>-\frac{N-2s}{2}
 \end{equation}
for every $r\in(0,\tilde R)$.
\end{Lemma}
\begin{pf}
From Corollary \ref{c:1sobolev_boundary}, \eqref{eq:ipoh}, and \eqref{eq:ipof},
it follows that
\begin{align*}
  \int_{B_r^+} t^{1-2s}|\nabla w|^2&\,dt\,dx- \kappa_s
  \int_{B_r'}\bigg(
  \dfrac{\lambda}{|x|^{2s}}w^2+hw^2+f(x,w)w\bigg)\,dx +
  \bigg(\frac{N-2s}{2r}\bigg)\int_{S_r^+}t^{1-2s}w^2dS\\
  & \geq \bigg(\frac{\kappa_s(\Lambda_{N,s}-\lambda)}2-C_h
  \kappa_sr^\e \bigg)
  \int_{B_r'}\dfrac{w^2}{|x|^{2s}}\,dx\\
  &\quad+\bigg(\frac{\Lambda_{N,s}-\lambda}{2(1+\Lambda_{N,s}) \widetilde
    S_{N,s}}-C_f\kappa_s|B_r'|^{\frac{2^*(s)-p}{2^*(s)}}\|w\|_{L^{2^*(s)}(B_r')}^{p-2}
  \bigg) \bigg(\int_{B_r'}|w|^{2^*(s)}\,dx\bigg)^{\!\!\frac2{2^*(s)}}
\end{align*}
for every $r\in(0,R_0)$. Since $\lambda<\Lambda_{N,s}$, from the above estimate
it follows that we can choose $\tilde R\in (0,R_0)$ sufficiently small
such that estimate (\ref{eq:13}) holds for $r\in (0,\tilde R)$ for
some positive constant
$\overline{C} >0$.
The proof of \eqref{eq:13bis} can be performed in a similar way, using
Lemmas \ref{l:hardy_boundary2} and \ref{l:sobolev_boundary}.
   Estimate (\ref{eq:13}), together with
(\ref{D(r)}) and (\ref{H(r)}), yields
(\ref{Nbelow}).
\end{pf}

\begin{Lemma} \label{l:stima_nu2} Let $\tilde R$ be as in Lemma \ref{l:stimasotto} and
  $\nu_2$ as in (\ref{eq:nu2}). Then there exist a positive constant
  $C_1>0$  and a function $g\in L^1(0,\tilde R)$,
  $g\geq 0$ a.e. in $(0,\tilde R)$, such that
$$
|\nu_2(r)|\leq C_1\left[{\mathcal
      N}(r)+\frac{N-2s}{2}\right]\Big(r^{-1+\e}+r^{-1+\frac{2s(p_0-2^*(s))}{p_0}}+g(r)\Big)
$$
for a.e. $r\in (0,\tilde R)$ and
$$
\int_0^r g(\rho)\,d\rho\leq
\frac{1}{1-\alpha}\|w\|_{L^{p}(B'_{\tilde R})}^{p(1-\alpha)}
  \, r^{ N \big(\frac{\alpha
    2^*(s)-2}{2^*(s)}-\frac{p\alpha-2}{p_0}\big)}
$$
for all $r\in (0,\tilde R)$ with some $\alpha$ satisfying $\frac{2}{p}<\alpha<1$.
\end{Lemma}

\begin{pf}
From \eqref{eq:ipoh} and  (\ref{eq:13}) we deduce that
\begin{align*}
\left|
\int_{B_r'} (2sh(x)+\nabla h(x)\cdot
    x)|w|^2\,dx\right|&\leq
2C_h r^\e
\int_{B_r'} \frac{|w|^2}{|x|^{2s}}\,dx\\
&\leq  2C_h\overline{C}^{-1}\,
r^{\e+N-2s}\left[D(r)+{\textstyle{\frac{N-2s}{2}}}H(r)\right],
\end{align*}
and, therefore, for any $r\in (0,\tilde R)$, we have that
\begin{align}\label{B00}
  \left|
\frac{\int_{B_r'} (2sh(x)+\nabla h(x)\cdot
    x)|w|^2\,dx}{\int_{S_r^+} t^{1-2s}w^2\,dS}\right|&\leq 2C_h\overline{C}^{-1}
  \,r^{-1+\e} \frac{D(r)+\frac{N-2s}{2}H(r)}{H(r)} \\
&\notag=
  2C_h\overline{C}^{-1}\, r^{-1+\e}\left[{\mathcal N}(r)+
    \frac{N-2s}{2}\right].
\end{align}
By \eqref{eq:ipof}, H\"older's inequality, and (\ref{eq:13}), for some
constant ${\rm const\,}={\rm const\,}(N,s,C_f)>0$ depending on
$N,s,C_f$, and for all $r\in(0,\tilde R)$, there holds
\begin{align*}
  \bigg|& \int_{B_r'} \big((N-2s)f(x,w)w-2NF(x,w)
  -2\nabla_xF(x,w)\cdot x\big)\, dx\bigg|
  \\
  &\leq {\rm const\,}\int_{B_r'}(|w|^{2}+|w|^{2^*(s)})\,dx
  \\
  &\leq {\rm const\,}\bigg(
  \Big(\frac{\omega_{N-1}}N\Big)^{\! \frac{2s}{N}}r^{2s}+\|w\|_{L^{2^*(s)}(B_{\tilde
      R}')}^{2^*(s)-2}\bigg) \bigg(\int_{B_r'}|w|^{2^*(s)}dx\bigg)^{\!\!\frac2{2^*(s)}}\\
 & \leq \frac{\rm const\,}{\overline C}\bigg(
  \Big(\frac{\omega_{N-1}}N\Big)^{\! \frac{2s}{N}}r^{2s} +
\Big(\frac{\omega_{N-1}}N\Big)^{\! \frac{2s(p_0-2^*(s))}{Np_0}}
r^{ \frac{2s(p_0-2^*(s))}{p_0}} \|w\|_{L^{p_0}(B_{\tilde
      R}')}^{2^*(s)-2}\bigg)
r^{N-2s}\left[D(r)+{\textstyle{\frac{N-2s}{2}}}H(r)\right]
\end{align*}
and hence
\begin{multline}\label{eq:20}
  \left|
\frac{\int_{B_r'}
    \big((N-2s)f(x,w)w-2NF(x,w) -2\nabla_xF(x,w)\cdot
    x\big)\, dx} {\int_{S_r^+} t^{1-2s}w^2\,dS}\right|
    \\
\leq
\frac{\rm const\,}{\overline C}
\bigg(
  \Big(\frac{\omega_{N-1}}N\Big)^{\! \frac{2s}{N}}r^{\frac{2s2^*(s)}{p_0}} +
\Big(\frac{\omega_{N-1}}N\Big)^{\! \frac{2s(p_0-2^*(s))}{Np_0}}
\|w\|_{L^{p_0}(B_{\tilde
      R}')}^{2^*(s)-2}\bigg)
r^{-1+\frac{2s(p_0-2^*(s))}{p_0} }
\left[\mathcal N(r)+{\textstyle{\frac{N-2s}{2}}}\right] .
\end{multline}
Let us fix $\frac{2}{p}<\alpha<1$. By H\"older's inequality and
(\ref{eq:13}),
\begin{align}\label{eq:21}
&  \bigg(\int_{B_r'}|w|^{p}\,dx\bigg)^{\!\!\alpha}=
\bigg(\int_{B_r'}|w|^{p-\frac2\alpha}|w|^{\frac2\alpha}\,dx\bigg)^{\!\!\alpha}\\
&\notag\leq
\bigg(\int_{B_r'}|w|^{2^*(s)\frac{p\alpha-2}{2^*(s)\alpha-2}}\,dx\bigg)^{\!\!\frac{\alpha2^*(s)-2}{2^*(s)}}
\bigg(\int_{B_r'}|w|^{2^*(s)}\,dx\bigg)^{\!\!\frac2{2^*(s)}}\\
&  \notag\leq
\Big(\frac{\omega_{N-1}}N\Big)^{\! \frac{\alpha
    2^*(s)-2}{2^*(s)}-\frac{p\alpha-2}{p_0}}r^{\! N \big(\frac{\alpha
    2^*(s)-2}{2^*(s)}-\frac{p\alpha-2}{p_0}\big)}
 \|w\|_{L^{p_0}(B_{\tilde
      R}')}^{p\alpha-2}
\frac{
r^{N-2s}}{\overline C}\left[D(r)+{\textstyle{\frac{N-2s}{2}}}H(r)\right]\\
  &\notag=
  \overline{C}^{-1}\Big(\frac{\omega_{N-1}}N\Big)^{\!\!\frac{\beta}N}
  r^{-1+\beta} \left[{\mathcal
      N}(r)+\frac{N-2s}{2}\right]\bigg(\int_{S_r^+}
  t^{1-2s}w^2\,dS\bigg)
\end{align}
for all $r\in(0,\tilde R)$, where
$\beta= N \big(\frac{\alpha
    2^*(s)-2}{2^*(s)}-\frac{p\alpha-2}{p_0}\big)>0$. From
\eqref{eq:ipof}, (\ref{eq:21}), and (\ref{Nbelow}), there exists some ${\rm
  const\,}={\rm const\,}(N,s,C_f)>0$ depending on $N,s,C_f$ such
that, for all $r\in(0,\tilde R)$,
\begin{multline}\label{eq:27}
  \left|
\frac{r\int_{\partial
      B_r'}\big(2F(x,w)-f(x,w)w\big)\, dS'} {\int_{S_r^+}
    t^{1-2s}w^2\,dS} \right| \leq {\rm const\,}
  \frac{r\int_{\partial B_r'}|w|^{p} \, dS'}{\int_{S_r^+}
    t^{1-2s}w^2\,dS}\\
\leq \frac{\rm
    const\,}{\overline{C}
}\Big(\frac{\omega_{N-1}}N\Big)^{\!\!\frac{\beta}N}
\left[{\mathcal
        N}(r)+\frac{N-2s}{2}\right]\frac{r^{\beta}  \int_{\partial B_r'}|w|^{p} \,
    dS'}{\Big(\int_{B_r'}|w|^{p} \, dx\Big)^{\!\alpha}}.
\end{multline}
By a direct calculation, we have that
\begin{equation}\label{eq:22}
\frac{r^{\beta}  \int_{\partial B_r'}|w|^{p} \,
    dS'}{\Big(\int_{B_r'}|w|^{p} \, dx\Big)^{\!\alpha}}
    =\frac{1}{1-\alpha}\left[\frac{d}{dr}\bigg(r^{\beta}
  \bigg( \int_{B_r'}|w|^{p} \, dx\bigg)^{\!\!1-\alpha}\bigg) -
  \beta\,  r^{-1+\beta}
  \bigg(\int_{B_r'}|w|^{p} \, dx\bigg)^{\!\!1-\alpha}\right]
\end{equation}
in the distributional sense and for a.e. $r\in (0,\tilde R)$.
Since
$$
\lim_{r\to0^+}
r^{\beta}
  \bigg( \int_{B_r'}|w|^{p} \, dx\bigg)^{\!\!1-\alpha}=0
$$
we deduce that the function
$$
r\mapsto
\frac{d}{dr}\bigg(r^{\beta}
  \bigg( \int_{B_r'}|w|^{p} \, dx\bigg)^{\!\!1-\alpha} \bigg)
$$
is integrable over $(0,\tilde R)$. Being
$$
r^{-1+\beta}
  \bigg(\int_{B_r'}|w|^{p} \, dx\bigg)^{\!\!1-\alpha} =o(r^{-1+\beta})
$$
as $r\to 0^+$, we have that also the function
$$
r\mapsto
r^{-1+\beta}
  \bigg(\int_{B_r'}|w|^{p} \, dx\bigg)^{\!\!\frac{p-2}{p}}
$$
is integrable over $(0,\tilde R)$. Therefore, by (\ref{eq:22}), we deduce that
\begin{equation}\label{eq:28}
g(r):=\frac{r^{\beta}  \int_{\partial B_r'}|w|^{p} \,
    {dS'}}{\Big(\int_{B_r'}|w|^{p} \, dx\Big)^{\!\alpha}}
\in L^1(0,\tilde R)
\end{equation}
and
\begin{align}\label{eq:29}
  0\leq \int_0^r g(\rho)\,d\rho&\leq
  \frac{1}{1-\alpha}\|w\|_{L^{p}(B'_{r})}^{p(1-\alpha)}
  \, r^{\beta}
\end{align}
for all $r\in (0,\tilde R)$. Collecting (\ref{B00}), (\ref{eq:20}),
(\ref{eq:27}), (\ref{eq:28}), and (\ref{eq:29}), we obtain the stated
estimate.
\end{pf}

\begin{Lemma} \label{l:stima_N_sopra} Let $\tilde R$ be as in Lemma \ref{l:stimasotto}
  and ${\mathcal N}$ as in (\ref{N(r)}). Then there exist a positive
  constant $C_2>0$   such that
\begin{equation} \label{Nabove}
{\mathcal N}(r)\leq C_2
\end{equation}
for all $r\in (0,\tilde R)$.
\end{Lemma}
\begin{pf}
By Lemma \ref{mono} and   Lemma \ref{l:stima_nu2}, we
 obtain
\begin{equation}\label{eq:40}
\bigg({\mathcal N}+\frac {N-2s}2\bigg)'(r)\geq
\nu_2(r)\geq
 -C_1\left[{\mathcal
      N}(r)+\frac{N-2s}{2}\right]\Big(r^{-1+\e}+r^{-1+\frac{2s(p_0-2^*(s))}{p_0}}+g(r)\Big)
\end{equation}
for a.e. $r\in (0,\tilde R)$.
Integration over $(r,\tilde R)$ yields
\begin{equation*}
{\mathcal N}(r)\leq -\frac {N-2s}2+
\left({\mathcal
      N}(\tilde R)+\frac{N-2s}{2}\right)
\exp\left(
C_1\bigg(
\frac{\tilde R^\e}{\e}+ \frac{p_0\tilde
R^{\frac{2s(p_0-2^*(s))}{p_0}}}{2s(p_0-2^*(s)) }
+\int_0^{\tilde R} g(\rho)\,d\rho\bigg)\right)
\end{equation*}
for any $r\in (0,\tilde R)$, thus proving estimate (\ref{Nabove}).
\end{pf}

\begin{Lemma} \label{gamma}
The limit
$$
\gamma:=\lim_{r\rightarrow 0^+} {\mathcal N}(r)
$$
exists and is finite.
\end{Lemma}
\begin{pf}
  By Lemmas \ref{l:stima_nu2} and \ref{l:stima_N_sopra}, the function
  $\nu_2$ defined in (\ref{eq:nu2}) belongs to $L^1(0,\tilde R)$.
  Hence, by Lemma \ref{mono}, ${\mathcal N}'$
  is the sum of a nonnegative function and of a $L^1$-function on
  $(0,\tilde R)$.  Therefore
$$
{\mathcal N}(r)={\mathcal N}(\tilde R)-\int_r^{\tilde R} {\mathcal N}'(\rho)\, d\rho
$$
admits a limit as $r\rightarrow 0^+$ which is necessarily finite in view of
(\ref{Nbelow}) and (\ref{Nabove}).
\end{pf}

\noindent The function $H$ defined in \eqref{H(r)} can be estimated as follows.
\begin{Lemma}\label{l:uppb}
  Let
  $\gamma:=\lim_{r\rightarrow 0^+} {\mathcal N}(r)$ be as in Lemma
  \ref{gamma} and $\tilde R$ as in Lemma \ref{l:stimasotto}.  Then
  there exists a constant $K_1>0$ such that
\begin{equation} \label{1stest}
H(r)\leq K_1 r^{2\gamma}  \quad \text{for all } r\in (0,\tilde R).
\end{equation}
Moreover, for any $\sigma>0$ there exists a constant
$K_2(\sigma)>0$ depending on $\sigma$ such that
\begin{equation} \label{2ndest} H(r)\geq K_2(\sigma)\,
  r^{2\gamma+\sigma} \quad \text{for all } r\in (0,\tilde R).
\end{equation}
\end{Lemma}

\begin{pf}
  By Lemma \ref{gamma}, ${\mathcal N}'\in L^1(0,\tilde r )$ and, by
  Lemma \ref{l:stima_N_sopra}, ${\mathcal N}$ is bounded, then from
  (\ref{eq:40}) and (\ref{eq:29}) it follows that
  \begin{equation} \label{qsopra}
{\mathcal N}(r)-\gamma=\int_0^r
    {\mathcal N}'(\rho) \, d\rho\geq -C_3 r^\delta
\end{equation}
for some constant $C_3>0$ and all $r\in(0, \tilde R)$, where
\begin{equation}\label{eq:45}
\delta=\min\bigg\{\e, \frac{2s(p_0-2^*(s))}{p_0},
 N \bigg(\frac{\alpha
    2^*(s)-2}{2^*(s)}-\frac{p\alpha-2}{p_0}\bigg)
\bigg\}>0.
\end{equation}
Therefore by (\ref{H'2}), \eqref{N(r)},
and (\ref{qsopra}) we deduce that, for all $r\in(0,\tilde R)$,
$$
\frac{H'(r)}{H(r)}=\frac{2\,{\mathcal N}(r)}{r}\geq
\frac{2\gamma}{r}-2C_3 r^{-1+\delta},
$$
which, after integration over the interval $(r, \tilde R)$, yields
(\ref{1stest}).

Since $\gamma=\lim_{r\rightarrow 0^+}
{\mathcal N}(r)$, for any $\sigma>0$ there exists $r_\sigma>0$ such
that ${\mathcal N}(r)<\gamma+\sigma/2$ for any $r\in (0,r_\sigma)$ and
hence
$$
\frac{H'(r)}{H(r)}=\frac{2\,{\mathcal N}(r)}{r}<\frac{2\gamma+\sigma}{r}
\quad \text{for all } r\in (0,r_\sigma).
$$
Integrating over the interval $(r,r_\sigma)$ and by continuity of $H$
outside $0$, we obtain (\ref{2ndest}) for some constant $K_2(\sigma)$
depending on $\sigma$.
\end{pf}

\section{The blow-up argument  }\label{sec:blow-up}

\noindent The main result of this section, which also contains Theorem \ref{t:asym-frac},
is the following theorem.
\begin{Theorem} \label{t:asym}
Let $w$ satisfy  \eqref{eq:wH8}, with $s,\lambda,h,f$ as in assumptions
\eqref{eq:ipo1},
\eqref{eq:ipoh}, and \eqref{eq:ipof}. Then, letting
  ${\mathcal N}(r)$ as in (\ref{eq:31}), there
there exists $k_0\in \N$, $k_0\geq1$, such that
\begin{equation}\label{eq:35}
\lim_{r\to 0^+}{\mathcal
      N}(r)=
-\frac{N-2s}{2}+\sqrt{\bigg(\frac{N-2s}
    {2}\bigg)^{\!\!2}+\mu_{k_0}(\lambda)}.
\end{equation}
Furthermore, if $\gamma$ denotes the limit in (\ref{eq:35}), $m\geq
  1$ is the multiplicity of the eigenvalue
  $\mu_{j_0}(\lambda)=\mu_{j_0+1}(\lambda)=\cdots=\mu_{j_0+m-1}(\lambda)$
  and
$\{\psi_i\}_{i=j_0}^{j_0+m-1}$ ($j_0\leq k_0\leq j_0+m-1$) is
  an $L^2({\mathbb
  S}^{N}_+;\theta_1^{1-2s})$-orthonormal basis for the eigenspace of
problem \eqref{eq:4} associated to $\mu_{k_0}(\lambda)$, then
\begin{align*}
&\tau^{-\gamma}w(0,\tau x)\to
|x|^{\gamma}\sum_{i=j_0}^{j_0+m-1} \beta_i\psi_{i}\Big(0,\frac x{|x|}\Big)
\quad \text{in }
C^{1,\alpha}_{\rm
    loc}(B_1'\setminus\{0\})  \quad \text{as }\tau\to 0^+,\\
&\tau^{-\gamma}w(\tau\theta)\to
\sum_{i=j_0}^{j_0+m-1} \beta_i\psi_{i}(\theta)\quad \text{in }
C^{0,\alpha}({\mathbb S}^{N}_+)  \quad \text{as }\tau\to 0^+,\\
& \tau^{-\gamma}w(0,\tau\theta')\to
\sum_{i=j_0}^{j_0+m-1} \beta_i\psi_{i}(0,\theta')\quad \text{in }
C^{1,\alpha}({\mathbb S}^{N-1}) \quad \text{as }\tau\to 0^+,
\end{align*}
and
\begin{equation*}
\tau^{1-\gamma}\nabla_x
  w(0,\tau\theta')\to \sum_{i=j_0}^{j_0+m-1}
  \beta_i\Big(\gamma\psi_{i}(0,\theta')\theta'+\nabla_{\SN} \psi_i(0,\cdot)(\theta')\Big)
  \quad \text{in } C^{0,\alpha}({\mathbb S}^{N-1}) \quad \text{as
  }\tau\to 0^+,
\end{equation*}
for some $\alpha\in(0,1)$, where
\begin{align*}
  \beta_i&= R^{-\gamma}
\int_{{\mathbb S}^{N}_+}\theta_1^{1-2s}w(R\,\theta)
  \psi_i(\theta)\,dS(\theta)
\\
  &\quad
+\kappa_s \int_{{\mathbb
      S}^{N-1}}\left[
\int_{0}^R\frac{h(t\theta')w(0,t\theta')+f(t\theta',
w(0,t\theta'))}{2\gamma+N-2s}\bigg(t^{2s-\gamma-1}-\frac{t^{\gamma+N-1}}{R^{2\gamma+N-2s}}\bigg)\,dt\right]
 \psi_i(0,\theta')\,dS(\theta'),
\end{align*}
 for all $R>0$ such that $\overline{B_{R}'}=
\{x\in\R^N:|x|\leq R\}\subset\Omega$
and $(\beta_{j_0},\beta_{j_0+1},\dots,\beta_{j_0+m-1})\neq(0,0,\dots,0)$.
\end{Theorem}

To prove Theorem \ref{t:asym}, we start by determining the asymptotic
profile of blowing up renormalized solutions to \eqref{eq:wHextended}.

\begin{Lemma}\label{l:blowup}
 Let  $w$ as in Theorem \ref{t:asym}. Let $\gamma:=\lim_{r\rightarrow 0^+} {\mathcal
    N}(r)$ as in Lemma \ref{gamma}. Then
\begin{itemize}
\item[\rm (i)] there exists $k_0\in \N$, $k_0\geq1$, such that
  $\gamma=-\frac{N-2s}{2}+\sqrt{\big(\frac{N-2s}
    {2}\big)^{2}+\mu_{k_0}(\lambda)}$;
\item[\rm (ii)] for every sequence $\tau_n\to0^+$, there exist a subsequence
$\{\tau_{n_k}\}_{k\in\N}$ and an eigenfunction $\psi$ of problem
\eqref{eq:4} associated to the eigenvalue $\mu_{k_0}(\lambda)$ such that
$\|\psi\|_{L^2({\mathbb
  S}^{N}_+;\theta_1^{1-2s})}=1$ and
\[
\frac{w(\tau_{n_k}z)}{\sqrt{H(\tau_{n_k})}}\to
|z|^{\gamma}\psi\Big(\frac z{|z|}\Big)
\]
 strongly in $H^1(B_r^+;t^{1-2s})$ and in  $C^{0,\alpha}_{\rm
    loc}(\overline{B_r^+}\setminus\{0\})$  for some $\alpha\in (0,1)$
  and all $r\in(0,1)$
and
\[
\frac{w(0,\tau_{n_k}x)}{\sqrt{H(\tau_{n_k})}}\to
|x|^{\gamma}\psi\Big(0,\frac x{|x|}\Big)
\]
in  $C^{1,\alpha}_{\rm
    loc}(B_1'\setminus\{0\})$.
\end{itemize}
\end{Lemma}
\begin{pf}
Let us set
\begin{equation}\label{eq:wtau}
w^\tau(z)=\frac{w(\tau z)}{\sqrt{H(\tau)}}.
\end{equation}
We notice that $\int_{S_1^+}t^{1-2s}|w^{\tau}|^2dS=1$. Moreover, by
scaling and (\ref{Nabove}),
\begin{multline}\label{eq:8bis}
  \int_{B_{1}^+} t^{1-2s}|\nabla w^\tau(z)|^2 dz -
\kappa_s
\int_{B_1'}\bigg(
\dfrac{\lambda}{|x|^{2s}}|w^\tau|^2+\tau^{2s}h(\tau
x)|w^\tau |^2\\+
\frac{\tau^{2s}}{\sqrt{H(\tau)}}
f\Big(\tau x,\sqrt{H(\tau)} w^\tau\Big)w^\tau\bigg)\,dx
  ={\mathcal N}(\tau)\leq C_2
\end{multline}
for every $\tau\in(0,\tilde R)$, whereas, from \eqref{eq:13bis},
\begin{multline}\label{eq:16}
  \mathcal N(\tau)\geq \frac{\tau^{-N+2s}}{H(\tau)}\bigg(
-\bigg(\frac{N-2s}{2\tau}\bigg)\int_{S_\tau^+}t^{1-2s}w^2dS
+\overline{C}\int_{B_\tau^+}
t^{1-2s}|\nabla w|^2\,dt\,dx\bigg)\\
=-\frac{N-2s}{2}+\overline{C} \int_{B_{1}^+} t^{1-2s}|\nabla w^\tau(z)|^2 dz
\end{multline}
for every $\tau\in(0,\tilde R)$. From \eqref{eq:8bis} and
\eqref{eq:16} we deduce that
\begin{equation}\label{eq:19}
  \{w^\tau\}_{\tau\in(0,\tilde R)}\quad\text{is bounded in }H^1(B_1^+;t^{1-2s}).
\end{equation}
Therefore, for any given sequence $\tau_n\to 0^+$, there exists a
subsequence $\tau_{n_k}\to0^+$ such that $w^{\tau _{n_k}}\weakly
\widetilde w$ weakly in $H^1(B_1^+;t^{1-2s})$ for some $\widetilde
w\in H^1(B_1^+;t^{1-2s})$.
 Due to compactness of the trace embedding \eqref{eq:tracecompact}, we
 obtain that
$\int_{S_1^+}t^{1-2s}|\widetilde w|^2dS=1$. In particular $\widetilde w\not\equiv
0$.

For every small  $\tau\in (0,\tilde R)$, $w^\tau$ satisfies
\begin{equation} \label{eqlam}
\begin{cases}
    \dive(t^{1-2s}\nabla  w^\tau)=0,&\text{in }B_1^+,\\
-\lim_{t\to 0^+}t^{1-2s}\frac{\partial w^\tau}{\partial
  t}=\kappa_s
\Big(\frac{\lambda}{|x|^{2s}}w^\tau+\tau^{2s}h(\tau
x)w^\tau+\frac{\tau^{2s}}{\sqrt{H(\tau)}}f(\tau x,\sqrt{H(\tau)}w^\tau)\Big), &\text{on }B_1',
\end{cases}
\end{equation}
in a weak sense, i.e.
\begin{multline}\label{eq:8tau}
\int_{B_1^+}t^{1-2s}\nabla
w^\tau\cdot\nabla\widetilde\varphi\,dt\,dx\\
=
\kappa_s
\int_{B_1'}\bigg(
\dfrac{\lambda}{|x|^{2s}}w^\tau+\tau^{2s}h(\tau
x)w^\tau+\frac{\tau^{2s}}{\sqrt{H(\tau)}}f\Big(\tau
x,\sqrt{H(\tau)}w^\tau\Big)\bigg)\widetilde\varphi(0,x)\,dx
\end{multline}
for all $\widetilde\varphi \in
H^1(B_1^+;t^{1-2s})$ s.t. $ \widetilde\varphi=0$ on $\S^N_+$ and,  for such  $\widetilde\varphi$,
  by \eqref{eq:ipof} and H\"older's inequality,
\begin{multline}\label{eq:39}
  \frac{\tau^{2s}}{\sqrt{H(\tau)}}\bigg|\int_{B_1'}f\Big(\tau
  x,\sqrt{H(\tau)}w^\tau(0,x)\Big)\widetilde\varphi(0,x)\,dx\bigg|\\
  \leq C_f\tau^{2s} \int_{B_1'}|w(0,\tau
  x)|^{p-2}|w^\tau(0,x)||\widetilde\varphi(0,x)|\,dx
  \\
  \leq C_f\|\mathop{\rm Tr}\widetilde\varphi\|_{L^p(B_1')}
\|w^\tau\|_{L^p(B_1')}\|w\|_{L^p(B_\tau')}^{p-2}\tau^{\frac{(N-2s)(2^*(s)-p)}{p}}
=o(1)\quad\text{as }\tau\to 0^+
\end{multline}
and, by \eqref{eq:ipoh} and Lemma \ref{l:hardy_boundary2},
\begin{multline}\label{eq:41}
 \tau^{2s}\left|\int_{B_1'}h(\tau
x)w^\tau\widetilde\varphi(0,x)\,dx\right|
\\\leq \frac{C_h\tau^\e}{\kappa_s\Lambda_{N,s}}
  \bigg(
 \int_{B_{1}^+} t^{1-2s}|\nabla w^\tau(z)|^2 dz+\frac{N-2s}2 \bigg)^{\!\!1/2}\!
  \bigg(
 \int_{B_{1}^+} t^{1-2s}|\nabla \widetilde\varphi (z)|^2 dz\bigg)^{\!\!1/2}
  =o(1)\text{ as }\tau\to 0^+.
\end{multline}
From (\ref{eq:39}), (\ref{eq:41}), and weak convergence
 $w^{\tau _{n_k}}\weakly
\widetilde w$ in $H^1(B_1^+;t^{1-2s})$, we can pass to
the limit in \eqref{eqlam} along the sequence $\tau _{n_k}$
and obtain  that $\widetilde w$ weakly solves
\begin{equation}\label{eq:extended_limit}
\begin{cases}
    \dive(t^{1-2s}\nabla \widetilde w)=0,&\text{in }B_1^+,\\
-\lim_{t\to 0^+}t^{1-2s}\frac{\partial \widetilde w}{\partial
  t}=\kappa_s
\frac{\lambda}{|x|^{2s}}\widetilde w, &\text{on }B_1'.
\end{cases}
\end{equation}
From \eqref{eq:ipof}, letting $q=\frac{p_0}{p-2}>\frac N{2s}$
with $p_0$ as in Lemma \ref{lem-upCrit}, we have that
\begin{align*}
&  \left\|
\frac{ \tau^{2s}}{\sqrt{H({\tau })}} \frac{f\big({\tau }
    x,\sqrt{H({\tau })}\,w^{\tau}(x)\big)}{w^{\tau }(x)}
\right\|_{L^q(B_2')}
\leq C_f\tau^{2s-\frac Nq}\bigg(\int_{B_{2\tau}}|w(x)|^{p_0}dx\bigg)^{\!\!1/q}\bigg)=O(1)
\end{align*}
as $\tau\to 0^+$.  Therefore from Lemma \ref{lem:reg} part (i) there holds
\begin{equation}\label{eq:33}
  w^{\tau _{n_k}}\to
  \widetilde w
\quad\text{in }C^{0,\alpha}_{\rm loc}(\overline{B_r^+}\setminus\{0\}),
\end{equation}
while Lemma \ref{lem:reg} part (ii) and Remark \ref{rem:unif_reg} imply that
\begin{equation}\label{eq:17}
  \nabla_xw^{\tau _{n_k}}\to
  \nabla_x\widetilde w,\quad\text{and}\quad
  t^{1-2s}\frac{\partial w^{\tau _{n_k}}}{\partial t}\to
  t^{1-2s}\frac{\partial \widetilde w}{\partial t}
  \quad\text{in }C^{0,\alpha}_{\rm loc}(\overline{B_r^+}\setminus\{0\})
\end{equation}
for some $\alpha\in (0,1)$ and all $r\in(0,1)$.
Reasoning as in \eqref{eq:39}, \eqref{eq:41}, we can prove that
\begin{equation}\label{eq:39kk}
  \frac{\tau^{2s}}{\sqrt{H(\tau)}}\int_{B_1'}f\Big(\tau
  x,\sqrt{H(\tau)}w^\tau\Big) w^\tau \,dx
=o(1)\quad\text{as }\tau\to 0^+
\end{equation}
and
\begin{equation}\label{eq:41kk}
 \tau^{2s}\int_{B_1'}h(\tau
x)|w^\tau|^2\,dx=o(1)\text{ as }\tau\to 0^+.
\end{equation}
Multiplying  equation \eqref{eqlam} with $w^\tau$, integrating in
$B_r^+$, and using
\eqref{eq:17}, \eqref{eq:39kk}, \eqref{eq:41kk}, we easily obtain that $\|w^{\tau
  _{n_k}}\|_{H^1(B_r^+;t^{1-2s})}\to \|\widetilde
w\|_{H^1(B_r^+;t^{1-2s})}$ for all $r\in(0,1)$, and hence
\begin{equation}\label{strongH1}
w^{\tau _{n_k}}\to \widetilde w\quad\text{in }H^1(B_r^+;t^{1-2s})
\end{equation}
for any $r\in (0,1)$.

For any $r\in (0,1)$ and $k\in \N$, let us define
the functions
\begin{align*}
D_k(r)=\frac{1}{r^{N-2s}} \bigg[\int_{B_r^+}
&t^{1-2s}|\nabla w^{\tau_{n_k}}|^2\,dt\,dx-
\kappa_s
\int_{B_r'}\bigg(
\dfrac{\lambda}{|x|^{2s}}|w^{\tau_{n_k}}|^2+\tau_{n_k}^{2s}h(\tau_{n_k}
x)|w^{\tau_{n_k}}|^2\\&+
\frac{\tau_{n_k}^{2s}}{\sqrt{H(\tau_{n_k})}}f\Big(\tau_{n_k}
x,\sqrt{H(\tau_{n_k})}w^{\tau_{n_k}}\Big)w^{\tau_{n_k}}\,\bigg)dx\bigg]
\end{align*}
and
\begin{equation*}
H_k(r)=\frac{1}{r^{N+1-2s}}\int_{S_r^+}t^{1-2s}|w^{\tau_{n_k}}|^2 \, dS.
\end{equation*}
Direct calculations yield
\begin{equation}\label{NkNw}
{\mathcal
    N}_k(r):=\frac{D_k(r)}{H_k(r)}=\frac{D(\tau_{n_k}r)}{H(\tau_{n_k}r)}
  ={\mathcal N}(\tau_{n_k}r) \quad \text{for all } r\in (0,1).
\end{equation}
From (\ref{strongH1}), \eqref{eq:39kk}, and \eqref{eq:41kk}, it follows that,
for any fixed $r\in (0,1)$,
\begin{equation} \label{convDk}
 D_k(r)\to \widetilde D(r),
\end{equation}
where
\begin{equation} \label{Dw(r)}
  \widetilde D (r)=
\frac{1}{r^{N-2s}} \bigg[\int_{B_r^+}
t^{1-2s}|\nabla \widetilde w|^2\,dt\,dx-
\kappa_s
\int_{B_r'}
\dfrac{\lambda}{|x|^{2s}}\widetilde w^2\,dx\bigg]
\quad \text{for all } r\in (0,1).
\end{equation}
On the other hand, by compactness of the trace embedding \eqref{eq:tracecompact}, we also have
\begin{equation} \label{convHk}
H_k(r)\to \widetilde H(r) \quad \text{for any fixed } r\in (0,1),
\end{equation}
where
\begin{equation} \label{Hw(r)}
\widetilde  H(r)=\frac{1}{r^{N+1-2s}}\int_{S_r^+}t^{1-2s}\widetilde
w^2 \, dS.
\end{equation}
From (\ref{eq:14}) it follows that $\widetilde
D(r)>-\frac{N-2s}2\widetilde H(r)$ for all
$r\in(0,1)$.  Therefore, if, for some $r\in(0,1)$, $\widetilde H(r)=0$ then
$\widetilde D(r)>0$, and passing to the limit in (\ref{NkNw}) should give a
contradiction with Lemma \ref{gamma}. Hence $H_w(r)>0$ for all
$r\in(0,1)$ and the function
\[
\widetilde {\mathcal N}(r):=\frac{\widetilde  D(r)}{\widetilde  H(r)}
\]
is well defined for $r\in (0,1)$. From (\ref{NkNw}), (\ref{convDk}), (\ref{convHk}), and
Lemma \ref{gamma}, we deduce  that
\begin{equation}\label{Nw(r)}
\widetilde {\mathcal
    N}(r)=\lim_{k\to \infty} {\mathcal
    N}(\t_{n_k}r)=\gamma
\end{equation}
for all $r\in (0,1)$.
Therefore $\widetilde {\mathcal N}$ is  constant in $(0,1)$ and hence $\widetilde {\mathcal
  N}'(r)=0$ for any $r\in (0,1)$.  By \eqref{eq:extended_limit} and Lemma
\ref{mono} with $h\equiv 0$ and $f\equiv0$, we obtain
\begin{equation*}
    \left(\int_{S_r^+}
  t^{1-2s}\left|\frac{\partial \widetilde w}{\partial \nu}\right|^2 dS\right) \cdot
    \left(\int_{S_r^+}
  t^{1-2s}\widetilde w^2\,dS\right)-\left(
\int_{S_r^+}
  t^{1-2s}\widetilde w\frac{\partial \widetilde w}{\partial \nu}\, dS\right)^{\!\!2}=0
\quad
  \text{for all } r\in (0,1),
\end{equation*}
which  implies  that $\widetilde w$ and $\frac{\partial \widetilde  w}{\partial \nu}$
have the same direction as vectors in $L^2(S_r^+;t^{1-2s})$ and hence
there exists a function $\eta=\eta(r)$ such that
$\frac{\partial \widetilde w}{\partial \nu}(r,\theta)=\eta(r)
\widetilde w(r,\theta)$ for all $r\in(0,1)$ and $\theta\in {\mathbb S}^N_+$.
After integration we obtain
\begin{equation} \label{separate}
\widetilde w(r,\theta)=e^{\int_1^r \eta(s)ds} \widetilde w(1,\theta)
=\varphi(r) \psi(\theta), \quad  r\in(0,1), \ \theta\in {\mathbb S}^N_+,
\end{equation}
where $\varphi(r)=e^{\int_1^r \eta(s)ds}$ and $\psi(\theta)=\widetilde w(1,\theta)$.
From \eqref{eq:extended_limit},  (\ref{separate}), and Lemma
\ref{l:separation},  it follows that, weakly,
$$
\begin{cases}
\frac1{r^N}\big(r^{N+1-2s}\varphi'\big)'\theta_1^{1-2s}\psi(\theta)+
r^{-1-2s}\varphi(r)\dive\nolimits_{{\mathbb
    S}^{N}}(\theta_1^{1-2s}\nabla_{{\mathbb S}^{N}}\psi(\theta))
=0,\\
-\lim_{\theta_1\to 0^+} \theta_1^{1-2s}\nabla_{{\mathbb
    S}^{N}}\psi(\theta)\cdot {\mathbf
  e}_1=\kappa_s\lambda\psi(0,\theta').
\end{cases}
$$
Taking $r$ fixed we deduce that $\psi$ is an eigenfunction of the
eigenvalue problem \eqref{eq:4}. If $\mu_{k_0}(\lambda)$ is the corresponding
eigenvalue then $\varphi(r)$ solves the equation
$$
\frac1{r^N}\big(r^{N+1-2s}\varphi'\big)'-\mu_{k_0}(\lambda)
r^{-1-2s}\varphi(r)=0
$$
i.e.
$$
\varphi''(r)+\frac{N+1-2s}r\varphi'-\frac{\mu_{k_0}(\lambda)}{r^2}\varphi(r)=0
$$
and hence $\varphi(r)$ is of the form
$$
\varphi(r)=c_1 r^{\sigma_{k_0}^+}+c_2 r^{\sigma_{k_0}^-}
$$
for some $c_1,c_2\in\R$, where
\begin{equation*}
  \sigma^+_{k_0}=-\frac{N-2s}{2}+\sqrt{\bigg(\frac{N-2s}
    {2}\bigg)^{\!\!2}+\mu_{k_0}(\lambda)}\quad\text{and}\quad
  \sigma^-_{k_0}=-\frac{N-2s}{2}-\sqrt{\bigg(\frac{N-2s}{2}
    \bigg)^{\!\!2}+\mu_{k_0}(\lambda
)}.
\end{equation*}
Since the function
$|x|^{\sigma_{k_0}^-}\psi(\frac{x}{|x|})\notin
L^2(B_1';|x|^{-2s})$ and hence $|z|^{\sigma_{k_0}^-}\psi(\frac{z}{|z|})
\notin  H^1(B_1^+;t^{1-2s})$  in virtue of Lemma
\ref{l:hardy_boundary2}, we deduce that $c_2=0$ and  $\varphi(r)=c_1
r^{\sigma_{k_0}^+}$. Moreover, from $\varphi(1)=1$, we obtain that $c_1=1$ and
then
\begin{equation} \label{expw}
\widetilde w(r,\theta)=r^{\sigma_{k_0}^+} \psi(\theta),  \quad
\text{for all }r\in (0,1)\text{ and }\theta\in {\mathbb S}^N_+.
\end{equation}
It remains to prove part (i).
From \eqref{expw} and the fact that $\int_{{\mathbb S}^N_+}\theta_1^{1-2s}\psi^2(\theta)dS=1$ it follows that
\begin{align*}
& \widetilde D (r)=
\frac{1}{r^{N-2s}} \bigg[\int_{B_r^+}
t^{1-2s}|\nabla \widetilde w|^2\,dt\,dx-
\kappa_s
\int_{B_r'}
\dfrac{\lambda}{|x|^{2s}}\widetilde w^2\,dx\bigg]  \\
&=r^{2s-N}(\sigma_{k_0}^+)^2\int_0^rt^{N-1-2s+2 \sigma_{k_0}^+}dt\\
&\quad+
r^{2s-N}\bigg(\int_0^rt^{N-1-2s+2 \sigma_{k_0}^+}dt\bigg)\bigg(
\int_{{\mathbb S}^{N}_+}\theta_1^{1-2s}|\nabla_{{\mathbb
        S}^{N}}\psi(\theta)|^2\,dS-\lambda
    \kappa_s \int_{{\mathbb
      S}^{N-1}}|\mathcal T\psi(\theta')|^2\,dS'
\bigg)\\
&=\frac{(\sigma_{k_0}^+)^2+\mu_{k_0}(\lambda)}{N-2s+2 \sigma_{k_0}^+}\,r^{2 \sigma_{k_0}^+}=
\sigma_{k_0}^+r^{2 \sigma_{k_0}^+}
\end{align*}
and
\begin{align*}
  \widetilde  H(r)=\int_{{\mathbb
      S}^{N}_+}\theta_1^{1-2s}\widetilde w^2(r\theta)\,dS=r^{2 \sigma_{k_0}^+},
\end{align*}
and hence from (\ref{Nw(r)}) it follows
that $\gamma=\widetilde{\mathcal N}(r)=\frac{\widetilde D(r)}{\widetilde H(r)}=\sigma_{k_0}^+$.
This completes the proof of the lemma.
\end{pf}

\noindent The following lemma describes the behavior of $H(r)$ as $r\to 0^+$.
\begin{Lemma} \label{l:limite}
Let $w$ satisfy \eqref{eq:wHextended}, $H$ be defined in
\eqref{H(r)}, and let $\gamma:=\lim_{r\rightarrow 0^+} {\mathcal
    N}(r)$ as in Lemma \ref{gamma}. Then the limit
\[
\lim_{r\to0^+}r^{-2\gamma}H(r)
\]
exists and it is finite.
\end{Lemma}
\begin{pf}
In view of (\ref{1stest}) it is sufficient to prove that the limit
exists. By (\ref{H(r)}), (\ref{H'2}), and Lemma~\ref{gamma} we have
$$
\frac{d}{dr} \frac{H(r)}{r^{2\gamma}} =-2\gamma r^{-2\gamma-1}
H(r)+r^{-2\gamma} H'(r) =2r^{-2\gamma-1} (D(r)-\gamma
H(r))=2r^{-2\gamma-1} H(r) \int_0^r {\mathcal N}'(\rho) d\rho.
$$
Integration over $(r,\tilde R)$ yields
\begin{equation}\label{inte}
  \frac{H(\tilde R)}{\tilde R^{2\gamma}}-
  \frac{H(r)}{r^{2\gamma}}=\int_r^{\tilde R} 2s^{-2\gamma-1}
  H(\rho) \left( \int_0^\rho \nu_1(t) dt \right) d\rho +\int_r^{\tilde R} 2\rho^{-2\gamma-1}
  H(\rho) \left( \int_0^\rho \nu_2(t) dt \right) d\rho
\end{equation}
where $\nu_1$ and $\nu_2$ are as in (\ref{eq:nu1}) and (\ref{eq:nu2}).
Since, by Schwarz's inequality, $\nu_1\geq 0$, we have that
$\lim_{r\to 0^+} \int_r^{\tilde R} 2 \rho^{-2\gamma-1} H(\rho) \left( \int_0^\rho
  \nu_1(t) dt \right) d\rho$
exists.  On the other hand, by
 (\ref{1stest}), Lemma \ref{l:stima_nu2}, and (\ref{Nabove}), we
deduce that
\begin{multline*}
  \left| \rho^{-2\gamma-1} H(\rho) \left( \int_0^\rho \nu_2(t) dt \right)
  \right|\leq K_1C_1\Big(C_2+\frac {N-2s}2\Big)\rho^{-1}\int_0^\rho\Big(
t^{-1+\e}+t^{-1+\frac{2s(p_0-2^*(s))}{p_0}}+g(t)\Big)\,dt\\
  \leq K_1C_1\Big(C_2+\frac{N-2s}2\Big)\rho^{-1}\bigg(
  \frac{\rho^\e}{\e}
+\frac{p_0 \rho^{\frac{2s(p_0-2^*(s))}{p_0}}}{2s(p_0-2^*(s)) }+
\frac{1}{1-\alpha}\|w\|_{L^{p}(B'_{\tilde R})}^{p(1-\alpha)}
  \, \rho^{ N \big(\frac{\alpha
    2^*(s)-2}{2^*(s)}-\frac{p\alpha-2}{p_0}\big)}
\bigg)
\end{multline*}
for all $\rho\in(0,\widetilde R)$, which proves that $\rho^{-2\gamma-1}
H(\rho) \left( \int_0^\rho \nu_2(t) dt \right)\in L^1(0,\widetilde R)$.
Hence both terms at the right hand side of
(\ref{inte}) admit a limit as $r\to 0^+$ thus completing the proof.
\end{pf}

\noindent From Lemma \ref{l:blowup}, the following pointwise estimate for
solutions to \eqref{eq:frac_eq} and  \eqref{eq:wH8} can be derived.
\begin{Lemma}\label{l:stima}
 Let  $w$
satisfying \eqref{eq:wH8}. Then there exist $C_4,C_5>0$ and
$\bar r\in(0,\tilde R)$ such that
\begin{enumerate}[\rm (i)]
\item $\sup_{S_{r}^+}|w|^2\leq
  \frac{C_4}{r^{N+1-2s}}\int_{S_r^+}t^{1-2s}|w(z)|^2\,dS$ for every $0<r<\bar r$,
\item $|w(z)|\leq C_5|z|^\gamma$ for all $z\in B_{\bar r}^+$ and in particular $|w(0,x)|\leq C_5|x|^\gamma$
 for all $x\in B_{\bar r}'$,
where $\gamma:=\lim_{r\rightarrow 0^+} {\mathcal
    N}(r)$ is as in Lemma \ref{gamma}.
\end{enumerate}
\end{Lemma}
\begin{pf}
We first notice that (ii) follows directly from (i) and
\eqref{1stest}. In order to prove (i), we argue by contradiction and
assume that there exists a sequence $\tau_n\to 0^+$ such that
\begin{equation*}
  \sup_{\theta\in{\mathbb S}^N_+}\Big|w\Big(\frac{\tau_n}2\theta\Big)\Big|^2>n H \Big (\frac{\tau_n}2 \Big)
\end{equation*}
with $H$ as in \eqref{H(r)}, i.e.
\begin{equation*}
  \sup_{x\in S_{1/2}^+}|w(\tau_nz)|^2>2^{N+1-2s}n\int_{S_{1/2}^+}t^{1-2s}w^2(\tau_nz)dS,
\end{equation*}
i.e., defining $w^\tau$ as in \eqref{eq:wtau}
\begin{equation}\label{eq:18}
  \sup_{x\in S_{1/2}^+}|w^{\tau_n}(z)|^2>2^{N+1-2s}n\int_{S_{1/2}^+}t^{1-2s}|w^{\tau_n}(z)|^2dS.
\end{equation}
From Lemma \ref{l:blowup}, along a subsequence $\tau_{n_k}$ we have
that $w^{\tau_{n_k}}\to
|z|^{\gamma}\psi\big(\frac z{|z|}\big)$ in  $C^{0,\alpha}_{\rm
    loc}(S_{1/2}^+)$, for some $\psi$ eigenfunction of problem
\eqref{eq:4}, hence passing to the limit in \eqref{eq:18} gives rise
to a contradiction.~\end{pf}

\noindent We will now prove that
$\lim_{r\to 0^+} r^{-2\gamma} H(r)$ is strictly positive.

\begin{Lemma} \label{l:limitepositivo} Under the same assumption as in
  Lemmas \ref{l:limite} and \ref{l:stima}, we have
\[
\lim_{r\to0^+}r^{-2\gamma}H(r)>0.
\]
\end{Lemma}
\begin{pf}
 For all $k\geq1$, let $\psi_k$ be as in \eqref{angular}, i.e. $\psi_k$
is a $L^2({\mathbb
  S}^{N}_+;\theta_1^{1-2s})$-normalized eigenfunction of problem
\eqref{eq:4} associated to the eigenvalue $\mu_k(\lambda)$ and
 $\{\psi_k\}_k$ is an orthonormal basis of $L^2({\mathbb
  S}^{N}_+;\theta_1^{1-2s}) $.
From Lemma \ref{l:blowup}
there exist
$j_0,m\in\N$, $j_0,m\geq 1$ such that $m$ is the multiplicity of the
eigenvalue $\mu_{j_0}(\lambda)=\mu_{j_0+1}(\lambda)=\cdots=\mu_{j_0+m-1}(\lambda)$ and
\begin{equation}\label{eq:30}
  \gamma=\lim_{r\rightarrow 0^+} {\mathcal N}(r)=-\frac{N-2s}{2}+\sqrt{\bigg(\frac{N-2s}
    {2}\bigg)^{\!\!2}+\mu_{i}(\lambda)},
  \quad i=j_0,\dots,j_0+m-1.
\end{equation}
Let us expand $w$ as
\begin{equation*}
w(z)=w(\tau \theta)=\sum_{k=1}^\infty\varphi_k(\tau)\psi_k(\theta)
\end{equation*}
where $\tau=|z|\in(0,R]$, $\theta=z/|z|\in{{\mathbb S}^{N}_+}$, and
\begin{equation}\label{eq:37}
  \varphi_k(\tau)=\int_{{\mathbb S}^{N}_+}\theta_1^{1-2s}w(\tau\,\theta)
  \psi_k(\theta)\,dS.
  \end{equation}
The Parseval identity yields
\begin{equation}\label{eq:17bis}
H(\tau)=\int_{{\mathbb
    S}^{N}_+}\theta_1^{1-2s}w^2(\tau\theta)\,dS=
\sum_{k=1}^{\infty}\varphi_k^2(\tau),\quad\text{for all }0<\tau\leq R.
\end{equation}
In particular, from
\eqref{1stest} and \eqref{eq:17bis} it follows that, for all $k\geq1$,
\begin{equation}\label{eq:23}
 \varphi_k(\tau)=O(\tau^{\gamma})\quad\text{as
 }\tau\to0^+.
\end{equation}
Equations \eqref{eq:wH8} and \eqref{angular} imply that, for every $k$,
\begin{equation*}
-\varphi_k''(\tau)-\frac{N+1-2s}{\tau}\varphi_k^\prime(\tau)+
\frac{\mu_k(\lambda)}{\tau^2}\varphi_k(\tau)=\zeta_k(\tau),\quad\text{in }(0,R),
\end{equation*}
where
\begin{equation}\label{eq:38}
  \zeta_k(\tau)=\frac{\kappa_s}{\tau^{2-2s}}\int_{{\mathbb
      S}^{N-1}}\big(h(\tau\theta')w(0,\tau\theta')+f(\tau\theta',
u(\tau\theta'))\big)  \psi_k(0,\theta')\,dS'.
\end{equation}
A direct calculation shows that, for some $c_1^k,c_2^k\in\R$,
\begin{equation}\label{eq:42}
\varphi_k(\tau)=\tau^{\sigma^+_k}
\bigg(c_1^k+\int_\tau^R\frac{t^{-\sigma^+_k+1}}{\sigma^+_k-\sigma^-_k}
\zeta_k(t)\,dt\bigg)+\tau^{\sigma^-_k}
\bigg(c_2^k+\int_\tau^R\frac{t^{-\sigma^-_k+1}}{\sigma^-_k-\sigma^+_k}
\zeta_k(t)\,dt\bigg),
\end{equation}
where
\begin{equation}\label{eq:46}
  \sigma^+_k=-\frac{N-2s}{2}+\sqrt{\bigg(\frac{N-2s}
    {2}\bigg)^{\!\!2}+\mu_k(\lambda)}\quad\text{and}\quad
  \sigma^-_k=-\frac{N-2s}{2}-\sqrt{\bigg(\frac{N-2s}{2}\bigg)^{\!\!2}+\mu_k(\lambda)}.
\end{equation}
From \eqref{eq:ipoh}, \eqref{eq:ipof},  Lemma \ref{l:stima},
\eqref{eq:30}, and the
fact that, in view of \eqref{firsteig_strict_in},
\[
2s+(p-2)\gamma=\frac{(N-2s)(2^*(s)-p)}{2}+(p-2) \sqrt{\bigg(\frac{N-2s}{2}\bigg)^{\!\!2}+\mu_{j_0}(\lambda)}>0,
\]
 we
deduce that, for all $i=j_0,\dots,j_0+m-1$,
\begin{equation}\label{eq:zeta}
  \zeta_{i}(\tau)=O(\tau^{-2+\tilde\delta+\sigma_{i}^+})\quad\text{as }\tau\to 0^+,
\end{equation}
with $\tilde\delta=\min\{\e,2s+(p-2)\gamma\}>0$. Consequently, the functions
\[
t\mapsto \frac{t^{-\sigma^+_{i}+1}}{\sigma^+_{i}-\sigma^-_{i}}
\zeta_{i}(t)\quad\text{and}\quad t\mapsto
\frac{t^{-\sigma^-_{i}+1}}{\sigma^-_{i}-\sigma^+_{i}} \zeta_{i}(t)
\]
belong to $L^1(0,R)$. Hence
\[
\tau^{\sigma^+_{i}}
\bigg(c_1^{i}+\int_\tau^R\frac{\rho^{-\sigma^+_{i}+1}}{\sigma^+_{i}-\sigma^-_{i}}
\zeta_{i}(\rho)\,d\rho\bigg)=o(\tau^{\sigma^-_{i}})\quad\text{as }\tau\to0^+,
\]
and then, by \eqref{eq:23}, there must be
\begin{equation*}
c_2^{i}=-\int_0^R\frac{t^{-\sigma^-_{i}+1}}{\sigma^-_{i}-\sigma^+_{i}}
\,\zeta_{i}(t)\,dt.
\end{equation*}
Using (\ref{eq:zeta}), we then deduce that
\begin{align}\label{eq:12}
  \tau^{\sigma^-_{i}}
  \bigg(c_2^{i}+\int_\tau^R\frac{t^{-\sigma^-_{i}+1}}{\sigma^-_{i}-\sigma^+_{i}}
  \zeta_{i}(t)\,dt\bigg)&=\tau^{\sigma^-_{i}}
  \bigg(\int_0^\tau
  \frac{t^{-\sigma^-_{i}+1}}{\sigma^+_{i}-\sigma^-_{i}}
  \zeta_{i}(t)\,dt\bigg)=O(\tau^{\sigma^+_{i}+\tilde\delta})
\end{align}
as $\tau\to0^+$.  From (\ref{eq:42}) and
(\ref{eq:12}), we obtain that,  for all $i=j_0,\dots,j_0+m-1$,
\begin{equation}\label{eq:24}
\varphi_i(\tau)=\tau^{\sigma^+_i}
\bigg(c_1^i+\int_\tau^R\frac{t^{-\sigma^+_i+1}}{\sigma^+_i-\sigma^-_i}
\zeta_i(t)\,dt+O(\tau^{\tilde\delta})\bigg)\quad\text{as }\tau\to0^+.
\end{equation}
Let us assume by contradiction that
$\lim_{\lambda\to0^+}\lambda^{-2\gamma}H(\lambda)=0$. Then, for all
$i\in\{j_0,\dots,j_0+m-1\}$,  (\ref{eq:30}) and (\ref{eq:17bis})
would imply that
\begin{equation*}
\lim_{\tau\to0^+}\tau^{-\sigma_{i}^+}\varphi_{i}(\tau)=0.
\end{equation*}
Hence, in view of \eqref{eq:24},
\[
c_1^{i}+\int_0^R\frac{t^{-\sigma^+_{i}+1}}{\sigma^+_{i}-\sigma^-_{i}}
\zeta_{i}(t)\,dt=0,
\]
which, together with (\ref{eq:zeta}), implies
\begin{align}\label{eq:13tris}
  \tau^{\sigma^+_{i}}
  \bigg(c_1^{i}+\int_\tau^R\frac{t^{-\sigma^+_{i}+1}}{\sigma^+_{i}-\sigma^-_{i}}
  \zeta_{i}(t)\,dt\bigg)= \tau^{\sigma^+_{i}}
  \int_0^\tau\frac{t^{-\sigma^+_{i}+1}}{\sigma^-_{i}-\sigma^+_{i}}
  \zeta_{i}(t)\,dt=O(\tau^{\sigma^+_{i}+\tilde\delta})
\end{align}
as $\tau\to0^+$. Collecting (\ref{eq:42}), (\ref{eq:12}), and
(\ref{eq:13tris}), we conclude that
\[
\varphi_{i}(\tau)=O(\tau^{\sigma^+_{i}+\tilde\delta})\quad\text{as
}\tau\to0^+\quad\text{for every }
i\in\{j_0,\dots,j_0+m-1\},
\]
namely,
\[
\sqrt{H(\tau)}\,(w^\tau,\psi)_{L^2({\mathbb
  S}^{N}_+;\theta_1^{1-2s})}=
O(\tau^{\gamma+\tilde\delta})\quad\text{as
}\tau\to0^+
\]
for every $\psi\in {\mathcal A}_0=\mathop{\rm span}\{
\psi_i\}_{i=j_0}^{j_0+m-1}$, where ${\mathcal A}_0$ is the eigenspace
of problem \eqref{eq:4} associated to the eigenvalue
  $\mu_{j_0}(\lambda)=\mu_{j_0+1}(\lambda)=\cdots=\mu_{j_0+m-1}(\lambda)$.
 From (\ref{2ndest}),
there exists $C(\tilde\delta)>0$ such that $\sqrt{H(\tau)}\geq
C(\tilde\delta)\tau^{\gamma+\frac\delta2}$ for $\tau$ small, and therefore
\begin{equation}\label{eq:26}
(w^\tau,\psi)_{L^2({\mathbb
  S}^{N}_+;\theta_1^{1-2s})}=
O(\tau^{\frac{\tilde\delta}2})\quad\text{as
}\tau\to0^+
\end{equation}
for every $\psi\in{\mathcal A}_0$.  From Lemma \ref{l:blowup},
 for every sequence $\tau_n\to0^+$, there exist a subsequence
$\{\tau_{n_k}\}_{k\in\N}$ and an eigenfunction $\widetilde
\psi\in{\mathcal A}_0$
\begin{equation}\label{eq:25}
  \int_{{\mathbb S}^{N}_+}\theta_1^{1-2s}\widetilde\psi^2(\theta)dS=1\quad\text{and} \quad
w^{\tau_{n_k}}\to \widetilde \psi\quad\text{in } L^2({\mathbb
  S}^{N}_+;\theta_1^{1-2s}).
\end{equation}
From (\ref{eq:26}) and (\ref{eq:25}), we infer that
\[
0=\lim_{k\to+\infty}(w^{\tau_{n_k}},\widetilde\psi)_{L^2({\mathbb S}^{N-1})}
=\|\widetilde\psi\|_{ L^2({\mathbb
  S}^{N}_+;\theta_1^{1-2s})}^2=1,
\]
thus reaching a contradiction.
\end{pf}

\noindent We can now completely describe the behavior of solutions to
\eqref{eq:wH8} near the
singularity, hence proving Theorem \ref{t:asym}.

\medskip\noindent \begin{pfn}{Theorem \ref{t:asym}} Identity
  (\ref{eq:35}) follows from part (i) of Lemma \ref{l:blowup}, thus
  there exists $k_0\in \N$, $k_0\geq 1$, such that $\gamma=\lim_{r\to
    0^+}{\mathcal N}(r)=-\frac{N-2s}{2}+\sqrt{\big(\frac{N-2s}
    {2}\big)^{2}+\mu_{k_0}(\lambda)}$.  Let us denote as $m$ the
  multiplicity of $\mu_{j_0}(\lambda)$ so that, for some $j_0\in\N$,
  $j_0\geq 1$, $j_0\leq k_0\leq j_0+m-1$,
  $\mu_{j_0}(\lambda)=\mu_{j_0+1}(\lambda)=\cdots=\mu_{j_0+m-1}(\lambda)$
and let $\{
\psi_i\}_{i=j_0}^{j_0+m-1}$ be an $L^2({\mathbb
  S}^{N}_+;\theta_1^{1-2s})$-orthonormal basis for the eigenspace
associated to $\mu_{k_0}(\lambda)$.

Let $\{\tau_n\}_{n\in\N}\subset (0,+\infty)$ such that
$\lim_{n\to+\infty}\tau_n=0$. Then, from part (ii) of Lemma
\ref{l:blowup} and Lemmas \ref{l:limite} and \ref{l:limitepositivo},
there exist a subsequence $\{\tau_{n_k}\}_{k\in\N}$ and $m$ real
numbers $\beta_{j_0},\dots,\beta_{j_0+m-1}\in\R$ such that
$(\beta_{j_0},\beta_{j_0+1},\dots,\beta_{j_0+m-1})\neq(0,0,\dots,0)$ and
\begin{align}\label{eq:32}
&\tau_{n_k}^{-\gamma}w(0,\tau_{n_k}x)\to
|x|^{\gamma}\sum_{i=j_0}^{j_0+m-1} \beta_i\psi_{i}\Big(0,\frac x{|x|}\Big)
\quad \text{in }
C^{1,\alpha}_{\rm
    loc}(B_1'\setminus\{0\})  \quad \text{as }k\to+\infty,\\
\label{eq:36}
&\tau_{n_k}^{-\gamma}w(\tau_{n_k}\theta)\to
\sum_{i=j_0}^{j_0+m-1} \beta_i\psi_{i}(\theta)\quad \text{in }
C^{0,\alpha}({\mathbb S}^{N}_+)  \quad \text{as }k\to+\infty,\\
\label{eq:47}&\tau_{n_k}^{-\gamma}w(0,\tau_{n_k}\theta')\to
\sum_{i=j_0}^{j_0+m-1} \beta_i\psi_{i}(0,\theta')\quad \text{in }
C^{1,\alpha}({\mathbb S}^{N-1}) \quad \text{as }k\to+\infty,
\end{align}
and
\begin{equation}\label{eq:48}
\tau_{n_k}^{1-\gamma}\nabla_x
  w(0,\tau_{n_k}\theta')\to\!\! \sum_{i=j_0}^{j_0+m-1}
  \beta_i(\gamma\psi_{i}(0,\theta')\theta'+\nabla_{\SN} \psi_i(0,\theta'))
   \quad\text{in } C^{0,\alpha}({\mathbb S}^{N-1}) \text{ as
  }k\to+\infty
\end{equation}
for some $\alpha\in(0,1)$.

 We now prove
that the $\beta_i$'s depend neither on the sequence
$\{\t_n\}_{n\in\N}$ nor on its subsequence
$\{\t_{n_k}\}_{k\in\N}$.

   Defining $\varphi_i$ and $\zeta_i$
as in (\ref{eq:37}) and \eqref{eq:38}, from
(\ref{eq:36}) it follows that, for any $i=j_0,\dots, j_0+m-1$,
\begin{equation}\label{eq:43}
\tau_{n_k}^{-\gamma}\varphi_i(\tau_{n_k}) =
\int_{{\mathbb S}^{N}_+}\theta_1^{1-2s}\frac{w(\tau_{n_k}\,\theta)}{\tau_{n_k}^{\gamma}}
  \psi_i(\theta)\,dS
\to\sum_{j=j_0}^{j_0+m-1} \beta_j
\int_{{\mathbb S}^{N}_+}\theta_1^{1-2s} \psi_j(\theta)
  \psi_i(\theta)\,dS=\beta_i
\end{equation}
as $k\to+\infty$.  As deduced in the proof of Lemma
\ref{l:limitepositivo}, for any $i=j_0,\dots, j_0+m-1$
and $\tau\in(0,R]$ there holds
\begin{align}\label{eq:44}
\varphi_i(\tau)&=
\tau^{\sigma^+_i}
\bigg(c_1^i+\int_\tau^R\frac{t^{-\sigma^+_i+1}}{\sigma^+_i-\sigma^-_i}
\zeta_i(t)\,dt\bigg)+
\tau^{\sigma^-_{i}}
  \bigg(\int_0^\tau
  \frac{t^{-\sigma^-_{i}+1}}{\sigma^+_{i}-\sigma^-_{i}}
  \zeta_{i}(t)\,dt\bigg)\\
\notag &=
\tau^{\sigma^+_i}
\bigg(c_1^i+\int_\tau^R\frac{t^{-\sigma^+_i+1}}{\sigma^+_i-\sigma^-_i}
\zeta_i(t)\,dt+O(\tau^{\tilde\delta})\bigg)\quad\text{as }\tau\to0^+,
\end{align}
for some $c_1^i\in\R$, where $\sigma^\pm_i$ are defined in \eqref{eq:46}.
Choosing $\tau=R$ in the first line of (\ref{eq:44}), we obtain
\[
c_1^i=R^{-\sigma^+_i}\varphi_i(R)-R^{\sigma^-_i-\sigma^+_i}\int_0^R
  \frac{s^{-\sigma^-_i+1}}{\sigma^+_i-\sigma^-_i}
  \zeta_i(s)\,ds.
\]
Hence (\ref{eq:44}) yields
\[
\tau^{-\gamma}\varphi_i(\tau)\to
R^{-\sigma^+_i}\varphi_i(R)-R^{\sigma^-_i-\sigma^+_i}\int_0^R
  \frac{t^{-\sigma^-_i+1}}{\sigma^+_i-\sigma^-_i}
  \zeta_i(t)\,dt+\int_0^R\frac{t^{-\sigma^+_i+1}}{\sigma^+_i-\sigma^-_i}
\zeta_i(t)\,dt\quad\text{as }\tau\to0^+,
\]
and therefore from (\ref{eq:43}) we deduce that
\begin{align*}
  \beta_i&= R^{-\gamma}
\int_{{\mathbb S}^{N}_+}\theta_1^{1-2s}w(R\,\theta)
  \psi_i(\theta)\,dS
\\
  &\quad-R^{-2\gamma-N+2s}\int_{0}^R\frac{\kappa_s \rho^{\gamma+N-1}}{2\gamma+N-2s}\bigg(
\int_{{\mathbb
      S}^{N-1}}\big(h(\rho\theta')w(0,\rho\theta')+f(\rho\theta',
w(0,\rho\theta'))\big)  \psi_i(0,\theta')\,dS'
\bigg) d\rho\\
  &\quad +\int_{0}^R\frac{\kappa_s \rho^{2s-\gamma-1}}{2\gamma+N-2s}\bigg(\int_{{\mathbb
      S}^{N-1}}\big(h(\rho\theta')w(0,\rho\theta')+f(\rho\theta',
w(0,\rho\theta'))\big)  \psi_i(0,\theta')\,dS'
 \bigg) d\rho .
\end{align*}
In particular the $\beta_i$'s depend neither on the sequence
$\{\tau_n\}_{n\in\N}$ nor on its subsequence
$\{\tau_{n_k}\}_{k\in\N}$, thus implying that the convergences in \eqref{eq:32},
(\ref{eq:36}), \eqref{eq:47}, and \eqref{eq:48} actually hold as $\tau\to 0^+$
and proving the theorem.~\end{pfn}

\section{Proof of Theorem \ref{t:asym-frac}}
Let $\mathcal D^{s,2}(\Omega)$ denote  the completion on
$C^\infty_{\rm c}(\Omega)$ with respect to the norm
$\|\cdot\|_{\mathcal D^{s,2}(\R^N)}$. Simple density arguments show that \eqref{eq:1} is equivalent to
\begin{equation}\label{eq:1pn}
(u,\varphi)_{\mathcal D^{s,2}(\R^N)}
=\int_\Omega \bigg(\frac{\lambda}{|x|^{2s}}u(x)+h(x)u(x)+f(x,u(x))\bigg)\varphi(x)\,dx,
\text{ for all }\varphi\in \calD^{s,2}(\Omega).
\end{equation}
Since $u\in \mathcal D^{s,2}(\R^N)$, we can let
$\mathcal H(u)\in
\mathcal D^{1,2}(\R^{N+1}_+;t^{1-2s})$ such that
$$
\int_{\R^{N+1}_+}t^{1-2s}\nabla \mathcal H(u)\cdot \nabla\phi\,dt\,dx=0,
\quad\text{for all }\phi \in C^\infty_{\rm c}({\R^{N+1}_+}),
$$
and $\mathcal H(u)
=u$ on $\R^N$ identified with $\de \R^ {N+1}_+$, i.e. $\mathcal H(u)$ weakly satisfies
$$
\begin{cases}
  \dive(t^{1-2s}\nabla \mathcal H(u))=0,&\text{in }\R^{N+1}_+,\\
\mathcal H(u)=u,&\text{on }\partial \R^{N+1}_+=\{0\}\times \R^N.
\end{cases}
$$
From \cite{CS} we have that
$$
\int_{\R^{N+1}_+}t^{1-2s}\nabla \mathcal
H(u)\cdot\nabla\widetilde\varphi\,dt\,dx=
\kappa_s(u, \widetilde\varphi)_{\mathcal
  D^{s,2}(\R^N)}\quad\text{for all }\widetilde\varphi\in  \mathcal D^{1,2}(\R^{N+1}_+;t^{1-2s}),
$$
where
$$
\kappa_s=\frac{\Gamma(1-s)}{2^{2s-1}\Gamma(s)},
$$
i.e.
$$
-\lim_{t\to 0^+}t^{1-2s}\frac{\partial \mathcal H(u)}{\partial
  t}=\kappa_s(-\Delta)^su(x),
$$
in a weak sense. Therefore $u\in \mathcal
  D^{s,2}(\R^N)$ weakly solves \eqref{eq:frac_eq} in $\Omega$ in the
  sense of \eqref{eq:1} if and only if its extension $w=\mathcal H(u)$
  satisfies
\begin{equation}\label{eq:extended}
\begin{cases}
    \dive(t^{1-2s}\nabla  w)=0,&\text{in }\R^{N+1}_+,\\
w =u,&\text{on }\R^{N},\\
-\lim_{t\to 0^+}t^{1-2s}\frac{\partial w}{\partial
  t}(t,x)=\kappa_s
\Big(\frac{\lambda}{|x|^{2s}}w+hw+f(x,w)\Big), &\text{on }\Omega,
\end{cases}
\end{equation}
in a weak sense, i.e. if  {for all }$\widetilde\varphi \in
\mathcal D^{1,2}({\R^{N+1}_+};t^{1-2s})$  such that $x\mapsto \widetilde\varphi(0,x)\in \mathcal D^{s,2}(\Omega)$ we have
\begin{equation}\label{eq:8}
 {\displaystyle{\int_{\R^{N+1}_+}t^{1-2s}\nabla w\cdot\nabla\widetilde\varphi\,dt\,dx=
\kappa_s
\int_\Omega\bigg(
\dfrac{\lambda}{|x|^{2s}}w+hw+f(x,w)\bigg)\widetilde\varphi\,dx}}.
\end{equation}
Since $\mathcal
D^{1,2}({\R^{N+1}_+};t^{1-2s})\hookrightarrow H^1(B_R^+;t^{1-2s})$ for all $R>0$, the result  follows from Theorem \ref{t:asym}.
\qed

\medskip\noindent \begin{pfn}{Theorem \ref{t:lambda0-asym-frac}}
The proof follows from the proof of Theorem \ref{t:asym-frac},
observing that, if $\lambda=0$, since no singularity occurs in the
coefficients of the equation, the convergences in \eqref{eq:33} and
\eqref{eq:17} hold in $C^{0,\alpha}(\overline{B_r^+})$ by
virtue of Lemma \ref{lem:reg}.
\end{pfn}

\medskip\noindent \begin{pfn}{Theorem \ref{t:sun}}
The proof follows directly from Theorem \ref{t:asym-frac}. Indeed, let
$u$ be a solution to \eqref{eq:frac_eq} satisfying
$u(x)=o(|x|^n)=o(1)|x|^n$ as $|x|\to 0$ for all $n\in
  \N$; if, by contradiction,  $u\not\equiv 0$ in $\Omega$,
  convergence \eqref{eq:limpsi} in
  Theorem \ref{t:asym-frac} would hold, thus contradicting the
  assumption
  that $u(x)=o(|x|^n)$ if $n>
-\frac{N-2s}{2}+\sqrt{ \left(\frac{2s-N}{2}  \right)^2
    +\mu_{k_0}(\l)  }$.
\end{pfn}

\section{Proof of Theorem \ref{t:unspm}}

Let $u\in \mathcal D^{s,2}(\R^N)$ be a weak
solution to \eqref{eq:frac_eq_unspm} in $\Omega$ such that $u\equiv 0$ on
a set $E\subset \Omega$ with $\mathcal L(E)>0$, where $\mathcal L$ denotes the
$N$-dimensional Lebesgue measure. By Lebesgue's density Theorem, for
a.e. point $x\in E$ there holds
\[
\lim_{r\to 0^+} \frac{\mathcal L (E\cap B(x,r))}{\mathcal L (B(x,r))}=1\quad\text{and}\quad
\lim_{r\to 0^+} \frac{\mathcal L ((\R^N\setminus E)\cap B(x,r))}{\mathcal L (B(x,r))}=0,
\]
where $B(x,r)=\{y\in\R^N:|y-x|<r\}$,
i.e. a.e. point of $E$ is a density point of $E$. Let $x_0$ be a
density point of $E$; then  for all $\e>0$ there exists $r_0=r_0(\e)\in(0,1)$
such that, for all $r\in(0,r_0)$,
\begin{equation}\label{eq:34}
\frac{\mathcal L ((\R^N\setminus E)\cap B(x_0,r))}{\mathcal L (B(x_0,r))}<\e.
\end{equation}
Assume by contradiction that $u\not\equiv 0$ in $\Omega$. Theorem
\ref{t:lambda0-asym-frac} implies that
\begin{equation}\label{eq:49}
r^{-\gamma}u(x_0+r(x-x_0))\to |x-x_0|^\gamma \psi
\bigg(0,\frac{x-x_0}{|x-x_0|}\bigg) \quad\textrm{as } r\to0^+,
\end{equation}
in  $C^{1,\alpha}(B_1(x_0,1))$, where
$\gamma=-\frac{N-2s}{2}+\sqrt{(\frac{N-2s}{2})^2+\mu_{k_0}}$ and
$\mu_{k_0}=\mu_{k_0}(0)\geq0$ is an eigenvalue of problem
\eqref{eq:4} with $\lambda=0$ and $\psi$ is a corresponding
eigenfunction. Since $u\equiv 0$ in $E$, by \eqref{eq:34} we have

\begin{align*}
\int_{B(x_0,r)}u^2(x)\,dx&=\int_{(\R^N\setminus E)\cap
  B(x_0,r)}u^2(x)\,dx\\
&\leq \bigg(\int_{(\R^N\setminus E)\cap
  B(x_0,r)}|u(x)|^{2^*(s)}dx\bigg)^{\!\!2/2^*(s)}|\mathcal L((\R^N\setminus E) \cap
B(x_0,r))|^{\frac{2^*(s)-2}{2^*(s)}}\\
&<\e^{\frac{2^*(s)-2}{2^*(s)}}|\mathcal
L(B(x_0,r))|^{\frac{2^*(s)-2}{2^*(s)}}
\bigg(\int_{(\R^N\setminus E)\cap
  B(x_0,r)}|u(x)|^{2^*(s)}dx\bigg)^{\!\!2/2^*(s)}
\end{align*}
for all $r\in(0,r_0)$, i.e. letting $u^r(x):=r^{-\gamma}u(x_0+r(x-x_0))$,
 \begin{align*}
\int_{B(x_0,1)}|u^r(x)|^2dx<\Big(\frac{\omega_{N-1}}N\Big)^{\!\frac{2^*(s)-2}{2^*(s)}}
\e^{\frac{2^*(s)-2}{2^*(s)}}\bigg(\int_{B(x_0,1)}|u^r(x)|^{2^*(s)}dx\bigg)^{\!\!2/2^*(s)}
\end{align*}
for all $r\in(0,r_0)$. Letting $r\to 0^+$, from \eqref{eq:49}  we have that
 \begin{align*}
\int_{B(x_0,1)}|x-x_0|^{2\gamma}&\psi^2\big(0,\tfrac{x-x_0}{|x-x_0|}\big)\,dx\\
&\leq
\Big(\frac{\omega_{N-1}}N\Big)^{\!\frac{2^*(s)-2}{2^*(s)}}
\e^{\frac{2^*(s)-2}{2^*(s)}}\bigg(\int_{B(x_0,1)}
|x-x_0|^{{2^*(s)}\gamma}\psi^{2^*(s)}\big(0,\tfrac{x-x_0}{|x-x_0|}\big) dx\bigg)^{\!\!2/2^*(s)}
\end{align*}
which yields a contradiction as $\e\to 0^+$.\qed

\medskip
 \textbf{Acknowledgements.} The authors would like to thank
 Prof. Enrico Valdinoci for his interest  and for taking
 their  attention to reference \cite{deFG} and the problem of unique
 continuation for sets of positive measure.


\begin{thebibliography}{99}

\bibitem{almgren} F. J. Jr. Almgren, {\it $Q$ valued functions
    minimizing Dirichlet's integral and the regularity of area
    minimizing rectifiable currents up to codimension two},
  Bull. Amer. Math. Soc., 8 (1983), no. 2, 327--328.


\bibitem{BCdP}
C. Brändle, E. Colorado, A. de Pablo, U. S\'anchez, {\it A concave-convex
    elliptic problem involving the fractional Laplacian},
  Proc. Roy. Soc. Edinburgh, 143A (2013), 39--71.


\bibitem{CaSi} X. Cabr\'e and Y. Sire, {\it Nonlinear equations for fractional Laplacians I: Regularity, maximum
    principles, and Hamiltonian estimates}, Annales de l'Institut
  Henri Poincaré (C) Non Linear Analysis, to appear, {\tt http://arxiv.org/abs/1012.0867}.


\bibitem{CS} L. Caffarelli, L. Silvestre, {\it An extension problem
    related to the fractional Laplacian},
 Comm. Partial Differential Equations 32 (2007), no. 7--9, 1245--1260.


\bibitem{carleman} T. Carleman,
{\it Sur un problème d'unicité pur les systèmes d'équations aux
  dérivées partielles à deux variables indépendantes}, Ark. Mat., Astr. Fys. 26 (1939), no. 17, 9 pp.

\bibitem{deFG} D. G. de Figueiredo, J.-P. Gossez, \emph{Strict monotonicity of
eigenvalues and unique continuation}, Comm. Partial Differential
Equations 17 (1992), no. 1--2, 339--346.

\bibitem{DNPV} E. Di Nezza, G. Palatucci, E. Valdinoci,
{\it Hitchhiker's guide to the fractional Sobolev spaces}, Bull. Sci. math. 136 (2012), no. 5, 521--573.

\bibitem{FKS} E. Fabes, C. Kenig, R. P. Serapioni,
{\it The local regularity of solutions of degenerate elliptic equations},
Comm. Partial Differential Equations 7 (1982), no. 1, 77--116.

\bibitem{fall-frac} M. M. Fall,  {\it
On a semilinear elliptic equation with fractional Laplacian and Hardy potential},
{\tt http://arxiv.org/abs/1109.5530}.

\bibitem{FW} M. M. Fall, T.  Weth, {\it Nonexistence results for a
    class of fractional elliptic boundary value problems}, J.  Funct.
  Anal., Vol. 263, no. 8 (2012) 2205--2227 .

\bibitem{FGL} E. Fabes, N. Garofalo, F.-H. Lin, {\it A partial answer
    to a conjecture of B. Simon concerning unique continuation},
  J. Funct. Anal. 88 (1990), no. 1, 194--210.

\bibitem{FFT} V. Felli, A. Ferrero, S. Terracini, {\it Asymptotic
    behavior of solutions to Schr\"odinger equations near an isolated
    singularity of the electromagnetic potential,}
 Journal of the European Mathematical Society 13 (2011), 119--174.

\bibitem{FFT2} V. Felli, A. Ferrero, S. Terracini, {\it On the
    behavior at collisions of solutions to Schr\"odinger equations with
    many-particle and cylindrical potentials,} Discrete Contin. Dynam. Systems. 32 (2012), 3895--3956.

\bibitem{FFT3} V. Felli, A. Ferrero, S. Terracini, {\it A note on
    local asymptotics of solutions to singular elliptic equations via
    monotonicity methods,} Milan
  J. Math 80 (2012), no. 1, 203--226.

\bibitem{FSV} A. Fiscella, R. Servadei, E. Valdinoci, {\it
    Asymptotically linear problems driven by the fractional}, Preprint 2012.

\bibitem{FLS} R. Frank, E. H. Lieb, R. Seiringer,
{\it Hardy-Lieb-Thirring inequalities for fractional Schr\"{o}dinger operators},
J. Amer. Math. Soc. 21 (2008), no. 4, 925-950.

\bibitem{GL} N. Garofalo, F.-H.  Lin, {\it Monotonicity properties of
    variational integrals, $A\sb p$ weights and unique continuation},
  Indiana Univ. Math. J.  35 (1986), no. 2, 245--268.

\bibitem{HLP} G. H. Hardy, J. E.  Littlewood, G.  P\'olya, {\it
    Inequalities}, 2d ed. Cambridge, at the University Press, 1952.

\bibitem{herbst} I. W. Herbst, {\it Spectral theory of the operator
    $(p^2+m^2)^{1/2}-Ze^2/r$},
Comm. Math. Phys. 53 (1977), no. 3, 285--294.

\bibitem{JK} D. Jerison, C. E. Kenig, {\it Unique continuation and absence of positive eigenvalues for Schrödinger operators}, Ann. of Math. (2) 121 (1985), no. 3, 463--494.

\bibitem{JLX} T. Jin, Y.Y. Li, J. Xiong, {\it On
a fractional Nirenberg problem, part I: blow up analysis and
compactness of solutions}, J. Eur. Math. Soc. (JEMS), to appear, {\tt arXiv:1111.1332v1}.


\bibitem{kurata} K. Kurata, {\it A unique continuation theorem for uniformly elliptic
 equations with strongly singular potentials}, Comm. Partial
 Differential Equations 18 (1993), no. 7-8, 1161--1189.

\bibitem{salsa} S. Salsa, {\it Partial differential equations in
    action. From modelling to theory},  Universitext. Springer--Verlag Italia, Milan, 2008.

\bibitem{silvestre} L. Silvestre, {\it Regularity of the obstacle
    problem for a fractional power of the Laplace operator},
Comm. Pure Appl. Math. 60 (2007), 67--112.

\bibitem{SV} Y. Sire, E. Valdinoci, {\it Fractional Laplacian phase
    transitions and boundary reactions: a geometric inequality and a
    symmetry result}, Journal of Functional Analysis 256 (2009), no. 6, 1842--1864.


\bibitem{wolff} T. H. Wolff, {\it A property of measures in $\R^ N$
    and an application to unique continuation},  Geom. Funct. Anal.,  2
    (1992), no. 2, 225--284.

\bibitem{yafaev} D. Yafaev, {\it
Sharp constants in the Hardy-Rellich inequalities},
J. Funct. Anal. 168 (1999), no. 1, 121--144.

\end{thebibliography}
\end{document}